\documentclass[twoside]{article}
\usepackage[a4paper]{geometry}
\usepackage[latin1]{inputenc} 
\usepackage[T1]{fontenc} 
\usepackage{RR}
\usepackage{hyperref}

\RRdate{January 2014}

\RRauthor{
Jean-Fr\'{e}d\'{e}ric Gerbeau
  \thanks[sfn]{Inria Paris-Rocquencourt, B.P. 105 Domain de Voluceau, 78153 Le Chesnay, France}%
	\thanks[sfnp]{UPMC Univ Paris 06, Laboratoire Jacques-Louis Lions, F-75005, Paris, France}
  \and
Damiano Lombardi
\thanksref{sfn}\thanksref{sfnp}
}

\authorhead{J-F. Gerbeau \& D. Lombardi}
\RRtitle{Paires de Lax approch\'ees pour l'int\'egration r\'eduite d'\'equations d'\'evolution non lin\'eaires}
\RRetitle{Approximated Lax Pairs for the Reduced Order Integration of Nonlinear Evolution Equations}
\titlehead{Approximated Lax Pairs for the Reduced Order Integration of Nonlinear Evolution Equations}
\RRresume{Un algorithme de r\'eduction de mod\`ele, appel\'e ALP, est propos\'e pour r\'esoudre de mani\`ere approch\'ee des \'equations d'\'evolution non lin\'eaires. Il est bas\'e sur une approximation de paires de Lax g\'en\'eralis\'ees.  Contrairement \`a d'autres m\'ethodes de r\'eduction de mod\`eles, comme la POD, la base sur laquelle la solution est cherch\'ee \'evolue selon une dynamique reli\'ee au probl\`eme. La m\'ethode est par cons\'equent bien adapt\'ee \`a des probl\`emes comportant des ondes progressives ou des propagations de fronts. Une autre diff\'erence avec d'autres m\'ethodes de r\'eduction de mod\`ele est qu'elle n'est pas bas\'ee sur une strat\'egie on-line / off-line. Nous montrons des exemples num\'eriques pour les \'equations de transport lin\'eaire, KdV et FKPP en dimension un et deux. 
}
\RRabstract{
A reduced-order model algorithm, called ALP, is proposed to solve nonlinear evolution partial differential equations. It is based on approximations of generalized Lax pairs. Contrary to other reduced-order methods, like Proper Orthogonal Decomposition, the basis on which the solution is searched for evolves in time according to a dynamics specific to the problem. It is therefore well-suited to solving problems with progressive front or wave propagation. Another difference with other reduced-order methods is that it is not based on an off-line / on-line strategy. Numerical examples are shown for the linear advection, KdV and FKPP equations, in one and two dimensions. 
 }
\RRkeyword{Reduced Order Modeling, Lax Pair, KdV, FKPP}
\RRmotcle{R\'eduction de mod\`ele, paires de Lax, KdV, FKPP}
\RRprojet{Reo}
\RCParis 

\usepackage{graphicx}
\usepackage{amsmath}
\usepackage{amssymb}
\usepackage{amsthm}
\usepackage{bbold}

\usepackage{color}
\usepackage[normalem]{ulem}
\definecolor{myorange}{rgb}{0.9568,0.4941,0.1961}
\definecolor{myred}{rgb}{0.9098,0.1294,0.2078}
\definecolor{myblue}{rgb}{0.0352,0.4981,0.6509}
\definecolor{mygreen}{rgb}{0.2235,0.6353,0.2588}
\usepackage{marvosym}

\newcommand{\R}{\ensuremath{\mathbb{R}}}
\newcommand{\Frac}[2]{\displaystyle{\frac{#1}{#2}}}
\newcommand{\Sum}[2]{\displaystyle{\sum_{#1}^{#2}}}
\theoremstyle{plain}
\newtheorem{prop}{Proposition}
\theoremstyle{definition}
\newtheorem*{rem}{Remark}

\begin{document}
\RRNo{8454}
\makeRR   
%

\section{Introduction}

This work is devoted to a method for solving time dependent nonlinear Partial Differential Equations (PDEs) by a reduced-order model (ROM). We are especially interested in problems exhibiting propagation phenomena. Those cases are well-known to be difficult to tackle with reduced order methods based on an off-line/on-line strategy like the Proper Orthogonal Decomposition (POD, see e.g. \cite{Sirovich_1}) or the Reduced Basis Method (RBM, see e.g. \cite{maday-ronquist-02,rozza-huynh-patera-07}). 

Unlike what is usually done, the method proposed in this work is based on a time dependent basis. In other words, the usual expansion $u(x,t) \approx \sum_{j=1}^{N_M} \beta_j(t) \phi_j(x)$
is replaced by 
$u(x,t) \approx \sum_{j=1}^{N_M} \beta_j(t) \phi_j(x,t)$.
Two questions have thus to be addressed: the definition of the basis and its propagation in time. 

To construct the basis, we propose to compute the eigenfunctions of a linear Schr\"odinger operator $\mathcal{L}(u_0) \cdot = -\Delta \cdot - \chi u_0 \cdot $ associated with the initial condition $u_0$ ($\chi$ being a given positive constant). This idea was inspired by recent works by Laleg, Cr\'epeau and Sorine~\cite{Sorine_1,laleg-08,laleg-crepeau-sorine-12} who proposed a signal processing technique, called  Semi-Classical Signal Analysis (SCSA). These authors showed in particular that these eigenfunctions could be used to obtain a parsimonious representation of the arterial blood pressure \cite{laleg-medigue-sorine-07}. 

Then we choose to propagate the basis in such a way that it remains an eigenbasis of the operator $\mathcal{L}(u(t)) \cdot = -\Delta \cdot - \chi u(t) \cdot $, where $u(t)$ is the solution of the PDE of interest at time $t$. The eigenfunctions satisfy a new evolution PDE associated with an operator $\mathcal{M}(u(t))$. In some particular cases, the operators $\mathcal{L}$ and $\mathcal{M}$ coincide with the ``Lax pairs'' introduced in \cite{lax-68}.
 
 Approximation with time dependent basis functions is of course not a new concept, see for example \cite{hughes1988space,petersen-farhat-tezaur-09,van2002space} among many others. The derivation of a system governing the evolution of an approximation basis has been recently introduced in the context of uncertainty quantification by \cite{Hou-Zhang-JCP1,Hou-Zhang-JCP2,Sapsis-Lermusiaux}. We are not aware of any other reduced order methods making use of time dependent basis. 
 
The structure of the work is as follows: 
in Section~\ref{sec:lax-pair}, some preliminary results and the links with Lax operators are briefly presented; Section~\ref{sec:alp-rom} is devoted to our reduced-order model algorithm, that will be called ALP, for Approximated Lax Pair. Some numerical tests are presented in Section~\ref{sec:num-exp} for the linear advection, Korteweg-de Vries and Fisher-Kolmogorov-Petrovski-Piskunov equations in one and two dimensions. In Section~\ref{sec:reduced-approx-scsa}, we show examples based on the Semi Classical Signal Analysis, which was used in a preliminary version of the present work~\cite{gerbeau:hal-00752810}.

\section{Preliminaries}
\label{sec:lax-pair}

\subsection{Time-dependent basis construction}

Let $\Omega$ be a bounded domain of $\R^d$ and $\chi$ a positive real number. Consider a real function $u(x,t)$, for $t\geq 0$ and $x=({x}_1,\dots,{x}_d)\in\Omega$. In the forthcoming sections, $u$ will be the solution of the PDE of interest. Throughout the paper, the functions will be assumed to have the regularity that justifies all the computations. 

Let $\mathcal{L}_{\chi}(u)$ be the Schr\"odinger operator associated with the potential $-\chi u$:
\begin{equation}
	\label{eq:def-Lu0-chi}
\mathcal{L}_{\chi}(u) \phi = - \Delta \phi - \chi u \phi,
\end{equation}
where $\Delta$ denotes the Laplacian in $d$ dimensions. For simplicity, the function $u(\cdot,t)$ will be denoted by $u(t)$ and the operator $\mathcal{L}_\chi(u(t))$ by $\mathcal{L}(t)$. The function $\phi$ is assumed to vanish on the boundary $\partial\Omega$. Other boundary conditions will be considered in the numerical tests. 

The operator $\mathcal{L}_{\chi}(u)$ is self-adjoint and the function $u$ is assumed to be regular enough so that $\mathcal{L}_{\chi}(u)$ has a continuous and compact inverse. A Hilbert basis of $L^2(\Omega)$ made of the eigenfunctions $(\phi_m(t))_{m>0}$ can therefore be defined as:
\begin{equation}
	\label{eq:eigenfunction-def}
	\mathcal{L}(t) \phi_m(t) = \lambda_m(t) \phi_m(t),
\end{equation}
where the (possibly negative) eigenvalues $\lambda_m(t)$ goes to $+\infty$ as $m\rightarrow \infty$.

Let $Q(t)$ be an orthogonal application ($Q^T Q = Q Q^T = Id$) such that 
$\phi_m(t) = Q(t) \phi_m(0), \forall m$. Taking the derivative with respect to $t$, we have:
\[
\partial_t \phi_m(t) = \partial_t Q(t) \phi_m(0) = \partial_t Q(t) Q^T(t)\phi_m(t).
\]
Thus, defining the operator $\mathcal{M}(t) = \partial_t Q(t) Q^T(t)$, the dynamics satisfied by the basis function is defined by:
\begin{equation}
\label{eq_dtPhi}
\partial_t \phi_m(t) = \mathcal{M}(t)\phi_m(t).
\end{equation}
Note that $\mathcal{M}^T = Q \partial_t Q^T = - \partial_t Q Q^T = - \mathcal{M}$, thus $\mathcal{M}$ is skew-symmetric.

To derive a relation between $\mathcal{L}$ and $\mathcal{M}$, we take the time derivative of $\mathcal{L}(t) \phi_m(t) = \lambda_m(t) \phi_m(t)$:
\[
\partial_t\mathcal{L} \phi_m + \mathcal{L}\mathcal{M} \phi_m = \partial_t \lambda_m \phi_m + \lambda_m \mathcal{M}\phi_m.
\]
Thus, defining the commutator $[\mathcal{L},\mathcal{M}] = \mathcal{L}\mathcal{M} - 
\mathcal{M}\mathcal{L}$, we obtain:
\begin{equation}
	\label{eq:lax-non-isospec}
(\partial_t \mathcal{L} + [\mathcal{L},\mathcal{M}]) \phi_m = \partial_t \lambda_m \phi_m	
\end{equation}
This equation will be instrumental for our algorithm. In particular, it will allow us to approximate operator $\cal{M}$ even when it is not known in closed-form.

\subsection{Links with the Lax Pairs}

Although this is not necessary for what follows, let us briefly show the links between $({\cal L},{\cal M})$ with the operators introduced by Lax in his seminal work~\cite{lax-68}. To integrate a class of nonlinear evolution PDEs, Lax introduced a pair of linear operators $\mathcal{L}(u)$ and $\mathcal{M}(u)$, where $u$ denotes the solution of the PDE.  These operators play the same role as in the previous section: the operator $\mathcal{L}(u)$ is defined as in \eqref{eq:def-Lu0-chi} and its eigenfunctions are propagated by $\mathcal{M}(u)$ as in \eqref{eq_dtPhi}. Lax focused on those particular cases when $\mathcal{L}(t)$ is orthogonally equivalent to $\mathcal{L}(0)$, i.e. when there exists $Q(t)$ orthogonal such that $\mathcal{L}(t) = Q(t) \mathcal{L}(0) Q^T(t)$. Then, defining as before $\mathcal{M}=\partial_t Q Q^T$, we have
\[
\partial_t Q^T \mathcal{L} Q + Q^T \partial_t \mathcal{L} Q  + Q^T \mathcal{L} \partial_t Q  = 0,
\]
left-multiplying by $Q$ and right-multiplying by $Q^T$, we obtain the \emph{Lax equation}:
\begin{equation}
	\label{eq:lax-isospec}
	\partial_t \mathcal{L} + [\mathcal{L},\mathcal{M}] = 0.
\end{equation}
A comparison of \eqref{eq:lax-non-isospec} and \eqref{eq:lax-isospec}, shows that in those cases the eigenvalues  satisfy $\partial_t \lambda_m = 0$. In other words, the eigenvalues $\lambda_m$ are ``first integrals of the motion''. When equation \eqref{eq:lax-isospec} is satisfied, operators $\mathcal{L}$ and $\mathcal{M}$ are said to be a \emph{Lax pair}. For some PDEs, it is possible to determine $\mathcal{M}$ in closed-form once $\mathcal{L}$ is chosen. A famous example is given by the Korteweg-de Vries equation (see below, Sections~\ref{sec:KdVNumerical} and \ref{sec:kdv-SCSA}, and equations \eqref{KdV_Operators-L}-\eqref{KdV_Operators-M}).

This formalism, which has close relations with the inverse scattering method, can be applied to a wide range of problems arising in many fields of physics (Camassa-Holm, Sine-Gordon, nonlinear Schr\"{o}dinger equations,\dots). In the huge literature devoted to Lax pairs, $\mathcal{M}(u)$ is generally used, or searched for, in closed-form (see for example \cite{Drazin_1} and the reference therein).  Most of the studies consider one dimensional domains and functions rapidly decreasing at infinity or periodic boundary conditions (see \cite{Fokas_1} for a theory on the finite line). 

In the present work, we are mainly interested in those cases when the eigenvalues are time dependent. Our work is therefore based on~\eqref{eq:lax-non-isospec} rather than on~\eqref{eq:lax-isospec}. In addition, we will not assume that operator $\mathcal{M}$ is explicitely known and we will consider bounded domains with Dirichlet or Neumann boundary conditions. Based on~\eqref{eq:lax-non-isospec}, we will propose an approximation of $\mathcal{M}$ and of the dynamics of $\lambda_m(t)$. For isospectral systems, \emph{i.e.} when \eqref{eq:lax-isospec} is satisfied, our method actually results in a numerical approximation of a Lax pair.  By abuse of language, we will keep the name ``Lax operators'' for $\mathcal{L}$ and $\mathcal{M}$ even for non-isospectral problems.

\section{Reduced-Order Modeling based on Approximated Lax Pairs (ALP)}
\label{sec:alp-rom}

We consider an evolution PDE set in $\Omega\times(0,T_{max})$:
\begin{equation}
	\label{DiffProb_1}
	\partial_{t} u = F(u),
\end{equation}
where $F(u)$ is an expression involving $u$ and its derivatives with respect to ${x}_1,\dots,{x}_d$. The problem is completed with an initial condition 
\begin{equation}
	\label{DiffProb_IC}
	u({x},0) = u_0({x}), \mbox{ for } {x}=({x}_1,\dots,{x}_d) \in \Omega.
\end{equation}
For the sake of simplicity, $u$ is assumed to vanish on the boundary $\partial\Omega$. Other boundary conditions will be considered in the numerical tests.

The solution $u({t})$ is searched for in a Hilbert space $V$ and approximated in $V_h$, a finite dimensional subspace of $V$, for example obtained by the finite element method (FEM). Let $(v_{j})_{j=1..N_{h}}$ denote a basis of $V_h$ and 
$\langle\cdot,\cdot\rangle$ the $L^2(\Omega)$ scalar product. 

\subsection{Reduced order approximation of the Lax operators}

The following proposition shows that it is possible to compute an approximation of $\mathcal{M}(u)$ in the space defined by the eigenfunctions of $\mathcal{L}_\chi(u)$ and to derive an evolution equation satisfied by the eigenvalues of $\mathcal{L}_\chi(u)$. 


\begin{prop}
	\label{prop:m-lambda}
	Let $u$ be a solution of equation \eqref{DiffProb_1}. Let $\mathcal{L}_\chi(u)$ be defined by \eqref{eq:def-Lu0-chi}. Let $N_M\in\mathbb{N}^\ast$. For $m\in \{1,\dots,N_M\}$, let $\lambda_m({t})$ be an eigenvalue of $\mathcal{L}_\chi(u({x},{t}))$, and $\phi_m({x},{t})$ an associated eigenfunction, normalized in $L^2(\Omega)$. Let $\mathcal{M}(u)$ be the operator defined in~\eqref{eq_dtPhi}.
Then the evolution of $\lambda_m$ is governed by
\begin{equation}
	\partial_{{t}} \lambda_m = - \chi \langle  F(u)\phi_m, \phi_m\rangle,
\label{eq:lambda-eq}
\end{equation}
and the evolution of $\phi_m$ satisfies, for $p\in\{1,\dots,N_M\}$,
\begin{equation}
	\label{eq:mode-eq}
	\langle \partial_{t} \phi_m, \phi_p \rangle  = M_{mp}(u),
\end{equation}
with
\begin{equation}
	\label{M_eq}
	\left\{
	\begin{array}{rcl}
M_{mp}(u) &=& \displaystyle{\frac{\chi}{\lambda_p - \lambda_m}} \langle  F(u)\phi_m, \phi_p\rangle, ~~\mbox{ if } p\neq m \mbox{ and } \lambda_p \neq \lambda_m,\\
M_{mp}(u) &=& 0, ~~ \mbox{ if } p=m  \mbox{ or } \lambda_p = \lambda_m.
\end{array}
\right.
\end{equation}
We will denote by $M(u) \in \R^{N_M\times N_M}$ the skew-symmetric matrix whose entries are defined by $M_{mp}(u)$.
\end{prop}

\textbf{Proof.}
Differentiating with respect to $t$ the equation satisfied by the $m$-th mode
\[
\mathcal{L}_\chi(u(x,t)) \phi_m(x,t) = \lambda_m(t) \phi_m(x,t),
\]
we get
\begin{equation}
\left(\mathcal{L}_\chi(u) - \lambda_m \mathcal{I}\right)\partial_{{t}} \phi_m = \partial_{{t}} \lambda_m \phi_m + \chi F(u)\phi_m.
\end{equation}

The scalar product is taken with a generic $\phi_p$, leading to:
\begin{equation}
\langle  \left(\mathcal{L}_\chi(u) - \lambda_m \mathcal{I}\right)\partial_{{t}} \phi_m,\phi_p\rangle  = \partial_{{t}} \lambda_m \langle  \phi_m,\phi_p\rangle  + \chi \langle  F(u)\phi_m, \phi_p\rangle.
\end{equation}
Using the self-adjointness of the operator and the orthonormality of the eigenfunctions, the following problem is obtained:
\begin{equation}
(\lambda_p-\lambda_m)\langle  \partial_{{t}} \phi_m, \phi_p\rangle  = \partial_{{t}} \lambda_m \delta_{mp} + \chi\langle  F(u)\phi_m, \phi_p\rangle.
\label{eq_Lax_1}
\end{equation}
Taking $p=m$, this proves \eqref{eq:lambda-eq}. In addition, the $L^2$ norm of $\phi_m$ being 1, $\langle \partial_{t} \phi_m, \phi_p \rangle = 0$, i.e. \eqref{M_eq}$_2$.

If $p\neq m$, but $\lambda_p = \lambda_m$ (multiple eigenvalues), we arbitrarilly set $M_{mp}(u) = 0$.
For $\lambda_p\neq \lambda_m$, we deduce from~\eqref{eq_Lax_1}:
\begin{equation}
\langle  \partial_{{t}} \phi_m, \phi_p\rangle  =  \frac{\chi }{\lambda_p - \lambda_m} \langle  F(u)\phi_m, \phi_p\rangle,
\label{dtMode}
\end{equation}
which completes the proof.
$\diamondsuit$

Equation \eqref{M_eq} gives an approximation of the operator $\mathcal{M}(u)$ on the basis defined by the modes at time ${t}$. This representation is convenient from a computational standpoint since it can easily be obtained from the expression $F(u)$ defining the PDE~\eqref{DiffProb_1}, without any \emph{a priori} knowledge of $\mathcal{M}(u)$. With this approximation of $\mathcal{M}(u)$, the evolution of the modes can be computed according to the nonlinear dynamics of the system. This is an important difference with standard reduced-order methods, like POD, where the modes are fixed once for all.

To set up a reduced order integration method, only a small number $N_M$ of modes will be retained. This number has to be chosen in order to represent the dynamics in a satisfactory way. A possible indicator of the quality of the approximation is given by the following quantity
\begin{equation}
\label{def-norm-dt-psi}
e(\phi_m({t}), N_M)  = \sum_{n=1}^{N_M} (M_{mn}(u({t})))^2,
\end{equation} 
which is an approximation of the $L^2$ norm of the time derivative of $\phi_m$:
\[
\int_{\Omega} \left(\partial_{{t}} \phi_m\right)^2\ d\Omega \approx \sum_{n,l=1}^{N_M} M_{ml}(u)M_{mn}(u) \langle  \phi_n,\phi_l\rangle  = \sum_{n=1}^{N_M} M_{mn}(u)^2 = e(\phi_m, N_M).
\]
By summing up over the modes, the Frobenius norm $\|\cdot \|_F$ of the representation of the evolution operator is recovered:
\begin{equation}
\label{frobDef}
\sum_{m=1}^{N_M} e(\phi_m({t}), N_M)  = \sum_{m,n=1}^{N_M} (M_{mn}(u({t})))^2  = 
\|M(u(t))\|_F^2.
\end{equation} 
This norm may be used as an error indicator for the dynamics recovery. This will be investigated in the numerical experiments presented in Section~\ref{sec:num-exp}.

\subsection{Reduced order approximation of the solutions}
\label{sec:reduced-approx}

The Hilbert basis defined by \eqref{eq:eigenfunction-def} is used to approximate the solution $u\in L^2(\Omega)$:
\begin{equation}
	\label{eq:u-modal}
	\tilde u({x}) = \sum_{m=1}^{N} \beta_m \phi_m(x),
\end{equation}

Another way of approximating the solution based on the mode squared and Deift-Trubowiz formula \cite{deift-trubowitz-79} was proposed and analyzed in~\cite{laleg-08,laleg-crepeau-sorine-12}. A preliminary version of the present work used this alternative representation~\cite{gerbeau:hal-00752810}. We now prefer using~\eqref{eq:u-modal} because of its generality. Nevertheless, the formula based on the mode squared deserves attention since it leads to a less expensive representation in some cases. The two approaches will be compared in~Section~\ref{sec:reduced-approx-scsa}.

\subsection{Reduced order dynamics}
\label{sec:rom-dynamics}

Proposition~\ref{prop:m-lambda} gives an approximated way to propagate the eigenmodes and the eigenvalues associated with a (generalized) Lax pair. Functions $u$ and $F(u)$ are approximated by \eqref{eq:u-modal} and 
\[
\tilde F(u) = \sum_{m=1}^{N_M} \gamma_m \phi_m,
\] 
respectively. Using these approximations in the PDE~\eqref{DiffProb_1}, the following holds: 
\[
\sum \dot\beta_m \phi_m + \beta_m \partial_t \phi_m = \sum \gamma_m \phi_m.
\]
Projecting this relation on $\phi_p$, and using~\eqref{eq:mode-eq}, the expression of the PDE on the reduced basis is obtained:
\[
	\dot\beta + M\beta = \gamma
\]
Defining $\Theta_{ij} = \langle \tilde F(u) \phi_j, \phi_i\rangle)$, 
\eqref{eq:lambda-eq} and \eqref{eq:mode-eq} are approximated by
\[
	\dot\lambda_i = - \chi \Theta_{ii}	,
\]
and, for $\lambda_i\neq \lambda_j$,
\[
	M_{ij} = \frac{\chi}{\lambda_j - \lambda_i} \Theta_{ij},
\]
respectively. 

\begin{rem}
It is also possible to derive the reduced order approximation of the Lax equation~\eqref{eq:lax-non-isospec}: $\frac{d\Lambda}{dt} + \chi \Theta = \Lambda M - M\Lambda$, from which the ordinary differential equations for $\beta$ and $\lambda_i$ can be straightforwardly deduced.
\end{rem}

The third order tensor $\langle \phi_k \phi_j, \phi_i \rangle$ is denoted by $T_{ijk}$. By definition:
\[
	\Theta_{ij} = \langle \tilde F(u) \phi_j, \phi_i\rangle = \sum_{k=1}^{N_M}\gamma_k T_{ijk}.	
\]
Computing the time derivative of $T_{ijk}$ gives:
\[
\dot T_{ijk} = \langle \partial_t \phi_k \phi_j,\phi_i \rangle
+ \langle \phi_k  \partial_t  \phi_j,\phi_i \rangle
+ \langle  \phi_k \phi_j, \partial_t \phi_i \rangle.
\]
Thus
\begin{equation}
	\dot T_{ijk} = \{M,T\} ^{(3)}_{ijk},	
\end{equation}
where 
\[
\{ M, T\}^{(3)}_{ijk} = \sum_{l=1}^{N_M}(  M_{li} T_{ljk} +  M_{lj} T_{ilk} +  M_{lk} T_{ijl}).
\]
For each specific equation, a relation linking $\gamma_i$ and $\beta=(\beta_j)$ will be also derived. For the time being, it is just generically denoted by $\gamma_i = \gamma_i(\beta)$.

For convenience, the expressions introduced in this section have been gathered in Table~\ref{table:notation}. To summarize, here is the set of equations which describes the dynamics in the reduced order space: 
\begin{equation}
\label{eq:rom-dynamics}
\left\{
	\begin{array}{rcl}
			\dot\beta_i + \Sum{m=1}{N_M} M_{im} \beta_m - \gamma_i & = & 0, \\
			\dot\lambda_i + \chi \Sum{m=1}{N_M} T_{iim} \gamma_m &=& 0, \\ 
		\dot T_{ijk} &=& \{M,T\}^{(3)}_{ijk}, \\ 
		M_{ij} &=& \Frac{\chi}{\lambda_j - \lambda_i} \Sum{m=1}{N_M} T_{ijm} \gamma_m, \\
		\gamma_i &=& \gamma_i(\beta), \\ 
	\end{array}
	\right.
\end{equation}
for $i,j,k = 1\dots N_M$. Relation~$\eqref{eq:rom-dynamics}_5$ will be made explicit in the examples given in Section~\ref{sec:num-exp}. 

As pointed out for example in~\cite{carlberg-bou-mosleh-farhat-11, ryckelynck-vincent-cantournet-12}, for any reduced order methods, it is generally to expensive to handle the nonlinearities of the equations by reconstructing the reduced order solution in the full-order space. 
Here, it is worth noticing that the integration is only done in the reduced-order space.

\begin{table}
	\begin{center}
\begin{tabular}{|c|c|}
\hline
Full Order Space & Reduced Order Space \\
\hline
$u$ & $\beta_m$ \\
$F(u)$ & $\gamma_m$ \\
$\mathcal{L}_\chi(u)$ & $\Lambda=\mbox{diag}(\lambda_m)$ \\
$\mathcal{M}(u)$ & $M_{mp}$ \\
$F(u)\cdot$ & $\Theta_{mp}$ \\
$(\partial_t \mathcal{L}_\chi + [\mathcal{L},\mathcal{M}]) \phi_m = \partial_t \lambda_m \phi_m$ & $\frac{d\Lambda}{dt} + \chi \Theta = \Lambda M - M\Lambda$ \\
$\partial_t u = F(u)$ & $\dot \beta + M \beta = \gamma$ \\
\hline	
\end{tabular}
\end{center}
\caption{Correspondence between the Full Order and Reduced Order spaces. The expression ``$F(u)\cdot$'' denotes the operator ``multiplication by the function $F(u)$''.}\label{table:notation}
\end{table}

\subsection{Numerical discretization of the reduced-order equation}
\label{sec:disc-alp-rom}
The ordinary differential equations system introduced in Section~\ref{sec:rom-dynamics} can be e.g. discretized by means of an implicit Runge-Kutta Gauss-Legendre method. For a generic vector $y(t)$ subject to the dynamics $\dot{y} = g(y)$, the method reads:
\begin{equation}
y^{(n+1)} = y^{(n)} + \delta t g\left( \frac{y^{(n)} + y^{(n+1)}  }{2}\right),
\label{eq:def-RK2}
\end{equation}
which is in general a nonlinear problem to be solved for $y^{(n+1)}$.

\subsection{The ALP algorithm}

\paragraph{Initialization} Let $u_0$ be the initial condition and let  $\epsilon_0>0$ be a prescribed tolerance. 
Compute a set of modes $(\phi^0_{m})_{m=1\dots N_M}$ and eigenvalues $(\lambda^0_{m})_{m=1\dots N_M}$ by solving
\[
\langle \nabla\phi^0_{m},\nabla v_{i}   \rangle - \chi \langle u_0 \phi^0_{m}, v_{i} \rangle = \lambda^0_{m} \langle \phi^0_{m},v_{i}\rangle, \mbox{ for } i=1,\dots,N_{h},	
\]
where $\chi$ is chosen such that $\| u_0 - \tilde u_0 \|_{L^2(\Omega)} \leq \epsilon_0$, with 
$\tilde u_0 = \sum_{j=1}^{N_M} \langle u_0, \phi_j^0 \rangle \phi_j^0$. The number of modes $N_M$ can be chosen in such a way that a criterion on the Frobenius norm of $M$ (see Eq.\eqref{eq::frobErr}) is satisfied.


\paragraph{Time evolution}

System~\eqref{eq:rom-dynamics} is discretized as follows:
\begin{equation}
	\label{eq:rom-dynamics-discret}
	\left\{
	\begin{array}{rcl}
\beta_i^{(n+1)} &=& \beta_i^{(n)} + \delta t \left( \gamma_i^{(n+1/2)} - \Sum{j=1}{N_M} M_{ij}^{(n+1/2) }\beta_j^{(n+1/2)} \right), \\
\mathcal{T}_{ijk}^{(n+1)} &=& \mathcal{T}_{ijk}^{(n)} + \delta t \left\lbrace \mathcal{T}^{(n+1/2)}, M^{(n+1/2)} \right\rbrace^{(3)}, \\
\lambda_i^{(n+1)} &=& \lambda_i^{(n)} -\chi \delta t \Sum{h=1}{N_M} \mathcal{T}_{iih}^{(n+1/2)} \gamma_h^{(n+1/2)},\\
M_{ji}^{(n+1/2)} &=& \Frac{\chi}{\lambda_j^{(n+1/2)}-\lambda_i^{(n+1/2)}} \Sum{h=1}{N_M}\mathcal{T}_{jih}^{(n+1/2)} \gamma_h^{(n+1/2)}, \\
 \gamma_i^{(n+1/2)} &=& \gamma_i(\beta^{(n+1/2)}).
\end{array}
\right . 
\end{equation}
where $\beta_i^{(n+1/2)} = (\beta_i^{(n)}+\beta_i^{(n+1)})/2$ and $T^{(n+1/2)} = (T^{(n+1)} + T^{(n)})/2$. 

Should other projection tensors be involved in the computation of the relation $\gamma(\beta)$, they would be updated as~$\mathcal{T}_{ijk}^{(n+1)}$ (see e.g. the case of the linear advection equation in Section~\ref{sec:num-exp}).

\subsection{Reduced order to full order transform}
\label{sec:reducedToFull}
In this section, the reconstruction of the solution in the full order space is addressed. This is done as a post-processing step, separated from the integration of the reduced-order model equations. 
Note that the reconstruction of the solution is more challenging than in classical ROM methods, since the basis evolves in time. 

Hereafter we present the simplest reconstruction method to compute the approximation of the solution in the full-order space, that is $\tilde{u}\approx \sum_{j=1}^{N_h} \hat{u}_j v_j$. We first note that equation~\eqref{eq:mode-eq} yields:
\begin{equation}
\label{eq:phi-with-residual}
\frac{\partial\phi_i}{\partial{t}} = \sum_{j=1}^{N_M} M_{ij}(u) \phi_j\ + r_i,
\end{equation}
where $r_i({t})\in\left[\mbox{span}(\phi_1({t}),\dots,\phi_{N_M}({t}))\right]^\perp$. Denoting by $B(t)$ the $N_h \times N_M$ matrix that represents the moving reduced order basis $(\phi_i(t))_{i=1..N_M}$ onto the fixed full order basis $(v_j)_{j=1..N_h}$:
\begin{equation}
\phi_j(x,t) = \sum_{i=1}^{N_h} B_{ij}(t) v_i(x).
\label{eq:repPhiOnV} 
\end{equation}

Equation~\eqref{eq:phi-with-residual} is approximated by neglecting the residual $r_i$ (the same notation is kept for simplicity):
\begin{equation}
	\label{eq:approx-dtB}
	\partial_t B = - B M.
\end{equation}
Various methods can be used to integrate this system. Here, we propose to use a simple two-step scheme.  First, \eqref{eq:approx-dtB} is integrated by means of a Crank-Nicolson scheme, then, to preserve orthonormality, a modified Gram-Schmidt algorithm is applied. The complexity of the Gram-Schmidt method is $2N_h\times N_M^2$, so that it is linear with respect to the full-order space dimension.
\begin{rem}
\label{rem:sub-spaceLim}
The reconstruction procedure proposed is straightforward and cheap from a computational standpoint, but it suffers from a limitation. Consider indeed \eqref{eq:approx-dtB}: it allows to take into account only the projection of the derivative of the modes on the modes themselves. When a linear update is performed, that means that the space spanned by the modes remains equal to the space they spanned at the previous time. It would be better to consider a correction of the form:
\begin{equation}
\partial_t B = - B M + W, \ \ W\in \R^{N_h\times N_M} \mbox{ such that } B^T G W = 0, 
\end{equation}
where $G$ is the Grammian matrix $[\langle v_j,v_i \rangle]$ and $W$ is determined by considering the evolution of the basis in the full-order space. This might be determined by considering that, from the integration of the reduced-order system, not only a representation of $u$ is available (namely $\beta$), but also of $\partial_t u$ (determined by $\gamma$). For the tests performed in the present work, the simple reconstruction proposed worked satisfactorily. The correction commented in this remark will be the object of further investigations.
\end{rem}






\section{Numerical Experiments}
\label{sec:num-exp} 

In this section, some numerical experiments are presented. The aim is to derive the reduced-order model for specific cases and to assess the numerical properties of the proposed algorithm. The first partial differential equation considered is a linear advection equation. It is a simple example of integrable system, i.e. for which the Lax pair is analytically known and satisfies equation~\eqref{eq:lax-isospec}. The second test case is performed on the Korteweg-de Vries equation that is a classical example of integrable system.

Then, the Fisher-Kolmogorov-Petrovski-Piskunov equation is addressed in one and two dimensions. This equation, which arises in many applications, features fronts propagation. Contrary to the linear advection and Korteweg-de Vries equations, it is not isospectral, i.e. \eqref{eq:lax-isospec} is not satisfied.

\subsection{Linear advection equation}
\label{sec:linWaveNumerical}

The ALP reduced-order method is first used to integrate the linear advection equation $\partial_t u + c \partial_x u = 0$. It is an integrable system: an exact Lax pair is given by the Schr\"{o}dinger operator and the constant operator $\mathcal{M} = - c \partial_x$. Following the ALP algorithm, this closed-form expression of $\mathcal{M}$ will not be used (see nevertheless the remark at the end of this section).

The function $u$ is approximated through an eigenfunction expansion of the form $u\approx \sum_{i=1} \beta_i(t) \phi_i(x,t)$, leading to:
\[
\sum_{j=1}^{N_M} \dot{\beta}_j \phi_j + \beta_j \partial_t \phi_j + c\sum_{j=1}^{N_M} \beta_j \partial_x \phi_j = 0,
\]
that becomes, after projection on a generic $\phi_i$:
\begin{equation}
	\label{eq:beta-advection}
\dot{\beta}_i + \sum_{j=1}^{N_M} M_{ij}\beta_j + c \sum_{j=1}^{N_M} D_{ij} \beta_j  = 0,
\end{equation}
where $D_{ij}:=\langle \partial_x \phi_j, \phi_i \rangle$ denotes the representation of the derivative operator in the reduced space. Identifying \eqref{eq:beta-advection} with $\eqref{eq:rom-dynamics}_1$ gives $\eqref{eq:rom-dynamics}_5$, i.e. the relation between $\beta$ and $\gamma$ specific to the advection equation:
\begin{equation}
	\label{eq:gamma-advection}
\gamma_i = -c \sum_{j=1}^{N_M} D_{ij}\beta_j.
\end{equation}
The evolution of matrix $D$ is governed by 
\begin{equation}
	\label{eq:D-dyn}
\dot D + [D,M] = 0.
\end{equation}
The system of equations to be solved in the reduced space is therefore $\eqref{eq:rom-dynamics}_{1-4}$, \eqref{eq:gamma-advection}, \eqref{eq:D-dyn}. It is solved with the numerical scheme described in Section~\ref{sec:disc-alp-rom}, with a time step $\delta t = 1/256$. The modes are computed with $\chi=150$. The initial condition is $u_0 = \exp(-250(x-0.25)^2)$, the advection velocity $c = 0.5$ and the final time is $T_{max} = 1$. 
\begin{figure}
\centerline{\hbox{\begin{tabular}{cc}
\includegraphics[height=7.5cm]{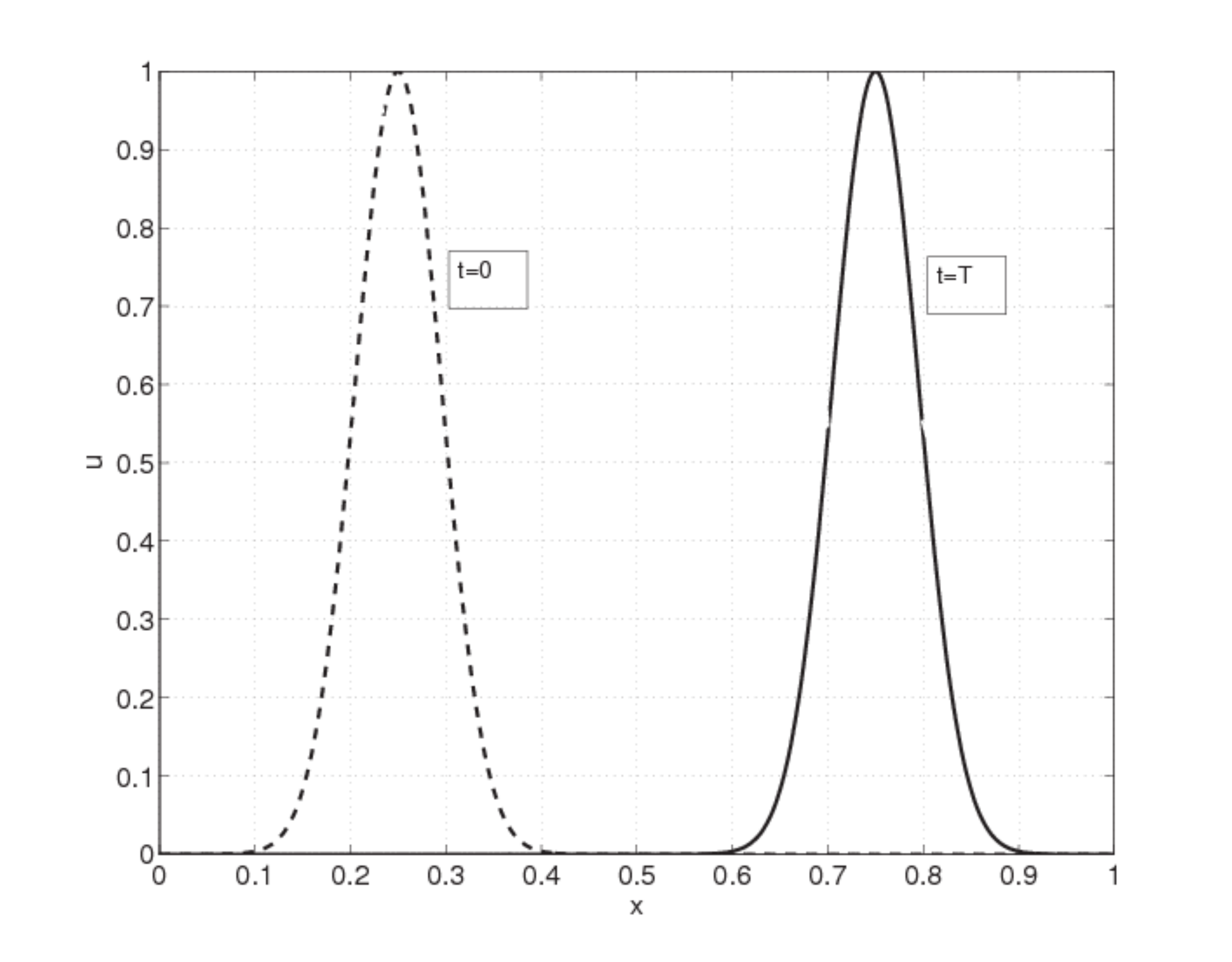} \\  (a) \\
\includegraphics[height=7.5cm]{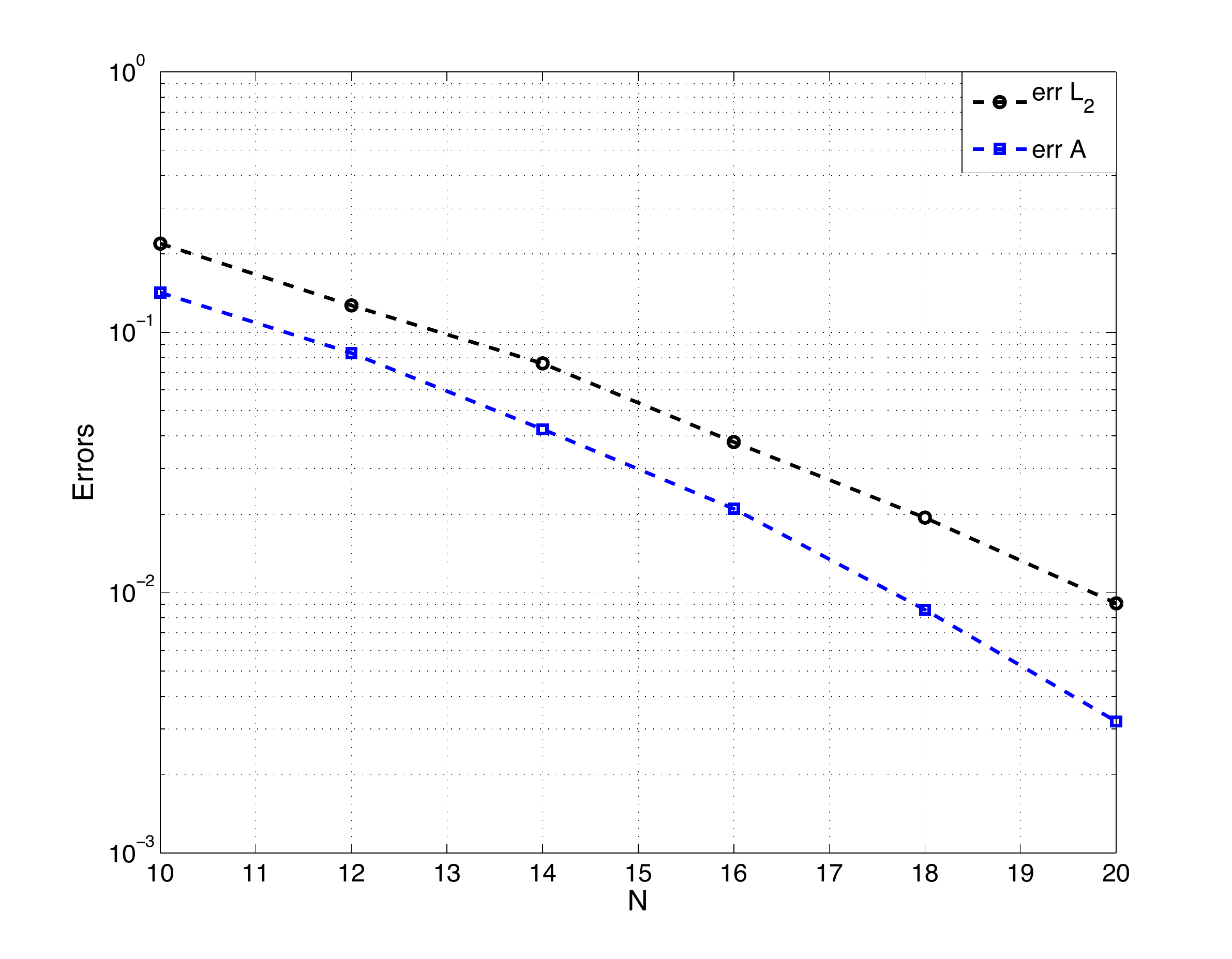} \\
(b) 
\end{tabular}}}
\caption{(a) Solutions at initial and final times, (b) errors in $L^2$ norm (black) and for the amplitude (blue) as a function of the number of modes in semi-logarithmic scale.}
\label{fig::LinWaveInitErr}
\end{figure}
In Fig.\ref{fig::LinWaveInitErr}.(a) the solution is represented at $t=0$ and $t=T_{max}$. The following error indicators are used to assess the quality of the solution:
\begin{eqnarray}
\varepsilon_{L2}^2 (t) := \frac{\int_{\Omega}(u - u_{ALP})^2\ d\Omega}{\int_{\Omega} u^2\ d\Omega}, \\
\varepsilon_{A} (t) := |\max(u) - \max(u_{ALP})|,
\end{eqnarray}
where $u_{ALP}$ is the reconstruction of the ROM solution, and $\varepsilon_{A} (t)$ assess the error in the peak amplitude. 

\begin{table}
\begin{center}
\begin{tabular}{cccc}
\hline
& $N_M$ & $\varepsilon_{L2}$ & $\varepsilon_{A}$ \\
\hline
& $10$ & $0.2193$ & $0.1421$ \\
& $12$ & $0.1267$ & $0.0832$ \\
& $14$ & $0.0759$ & $0.0424$ \\
& $16$ & $0.0379$ & $0.0210$ \\
& $18$ & $0.0194$ & $0.0086$ \\
& $20$ & $0.0091$ & $0.0032$ \\
\hline
\end{tabular}\caption{Error indicators for the linear advection test case as a function of the number of modes used to discretize the equations: first column $N_M$ is the number of modes, the second and the third ones the errors in $L^2$ norm and in the wave amplitude.}
\label{table::LinWaveErr}
\end{center}
\end{table}

These two error indicators were computed and evaluated by varying the number of modes used to discretize the equations in the reduced space. In Fig.\ref{fig::LinWaveInitErr}.(b) the errors $\overline{\varepsilon}_{L2} = \int_0^T \varepsilon_{L2} \ dt $ and $\| \varepsilon_A\|_{\infty}$ are plotted in a semi-logarithmic scale as a function of the number of modes used. Note the exponential convergence of the method and the fact that the error in the peak position is weakly dependent on the number of modes used. When $N_M=20$ modes, the method has roughly the same error as the Lax-Friedrichs scheme with twice as many iterations in time, optimal CFL and $1000$ space points.

\begin{figure}
\centerline{\hbox{\begin{tabular}{cc}
\includegraphics[width=9cm]{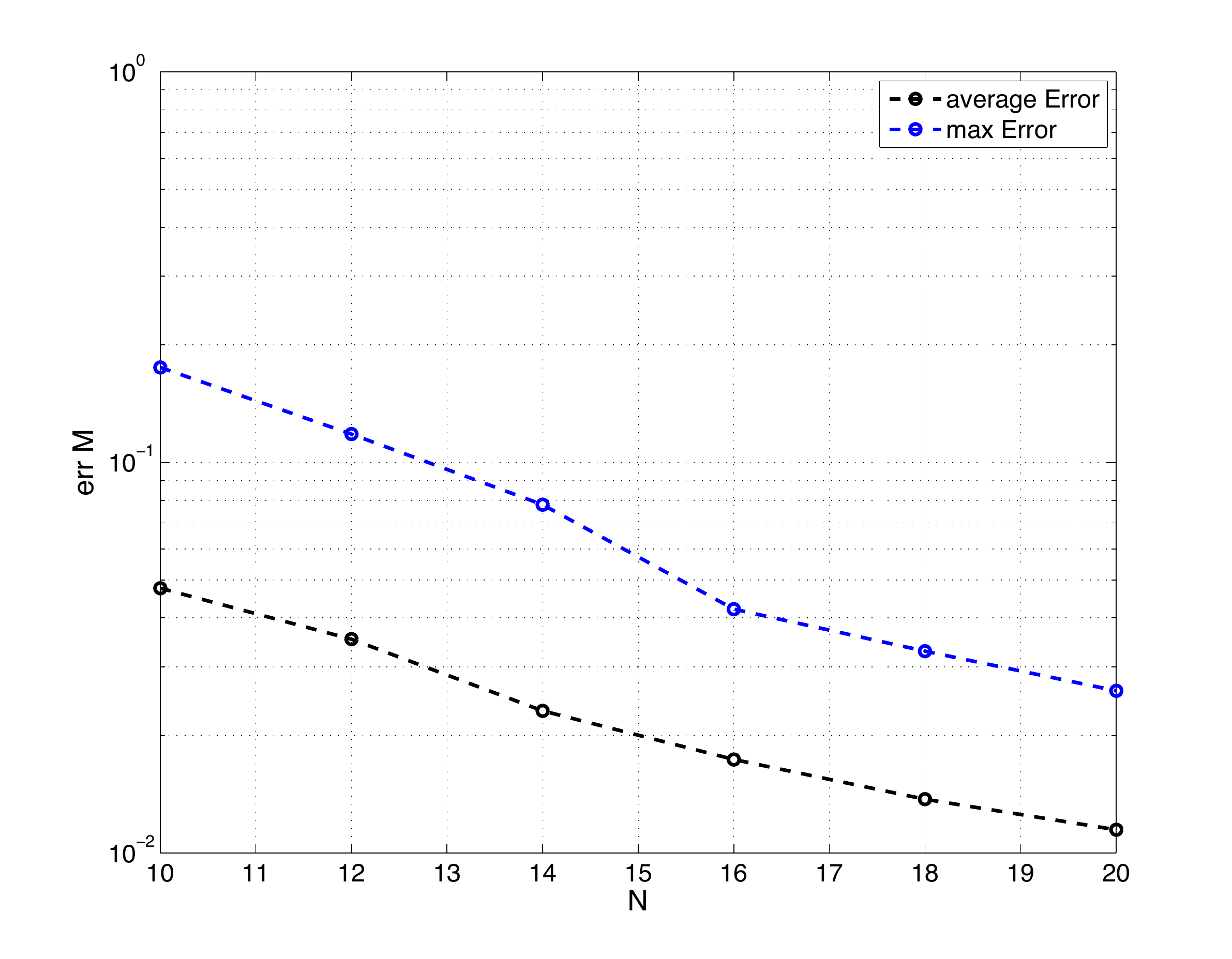} \\  (a) \\
\includegraphics[width=9cm]{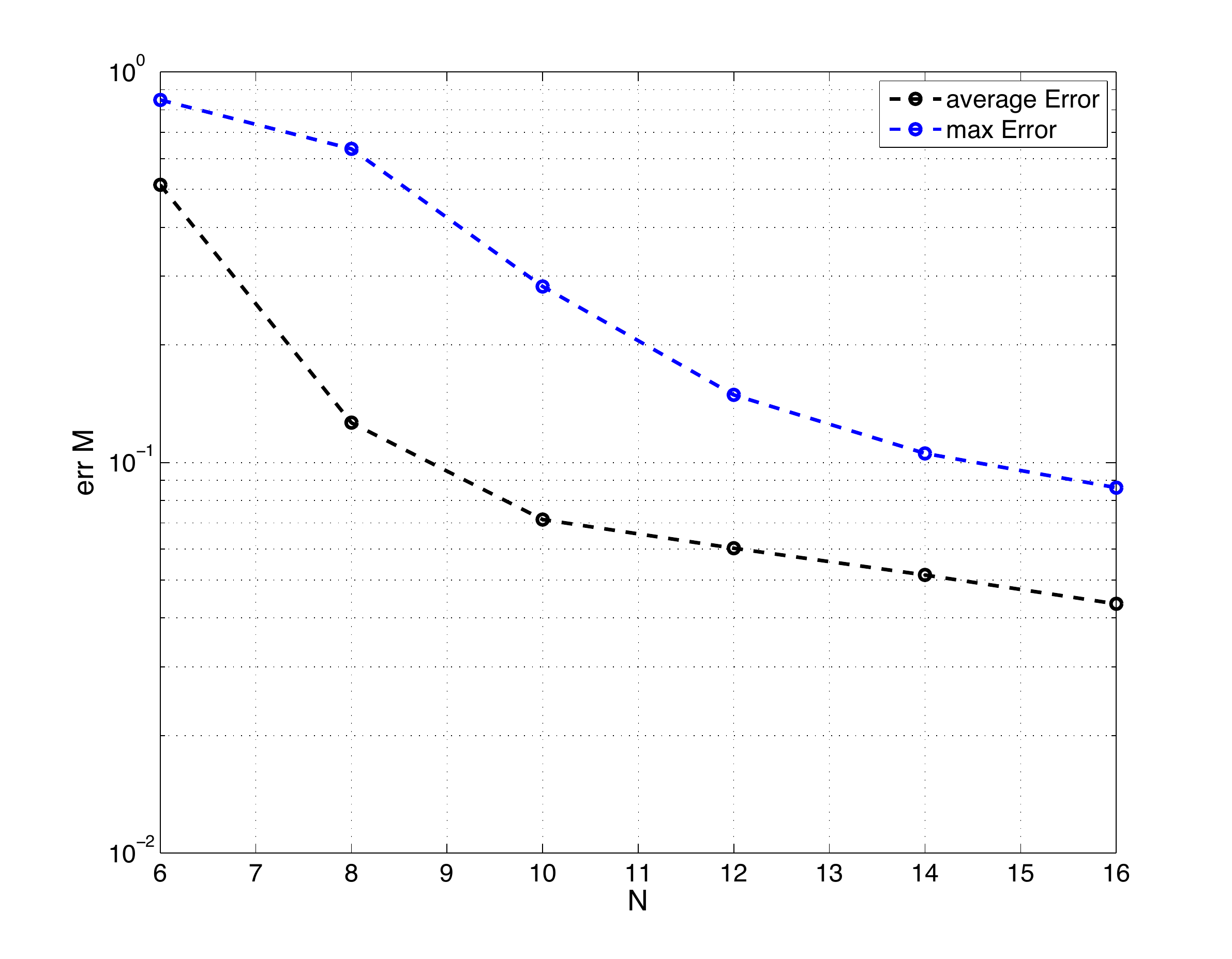} \\  (b) 
\end{tabular}}}
\caption{Time average and maximum of the Frobenius norm error indicator (see \eqref{eq::frobErr}) of $M$ as a function of the number of modes used in semi-logarithmic scale for (a) the linear advection equation (section~\ref{sec:linWaveNumerical}), (b) the 1D FKPP equation (section~\ref{sec:fkpp-1d}).}
\label{fig::Frobenius}
\end{figure}

The Frobenius norm of the matrix $M$ can be used as an intrinsic error indicator to evaluate the quality of the dynamical reconstruction of the solution (see \eqref{frobDef}). Let us define the error indicator:
\begin{equation}
\varepsilon_M(t,N_M) := \frac{| \| M_{N_M} \|_F - \|M_{\infty}\|_F}{\|M_{\infty}\|_F},
\label{eq::frobErr}  
\end{equation}
where $\|M_{\infty}\|_F$ is the Frobenius norm of the operator computed by means of $N_M=50$ modes.
In Fig.\ref{fig::Frobenius}.(a) the time average and the maximum of the Frobenius norm error indicator is shown as a function of the number of modes, in semi-logarithmic scale, for the linear advection equation test case. This plot suggests that the Frobenius norm criterion might be a good estimator to evaluate the convergence of the ROM towards the solution.

\begin{figure}
\centerline{\hbox{\begin{tabular}{cc}
\includegraphics[height=7.5cm]{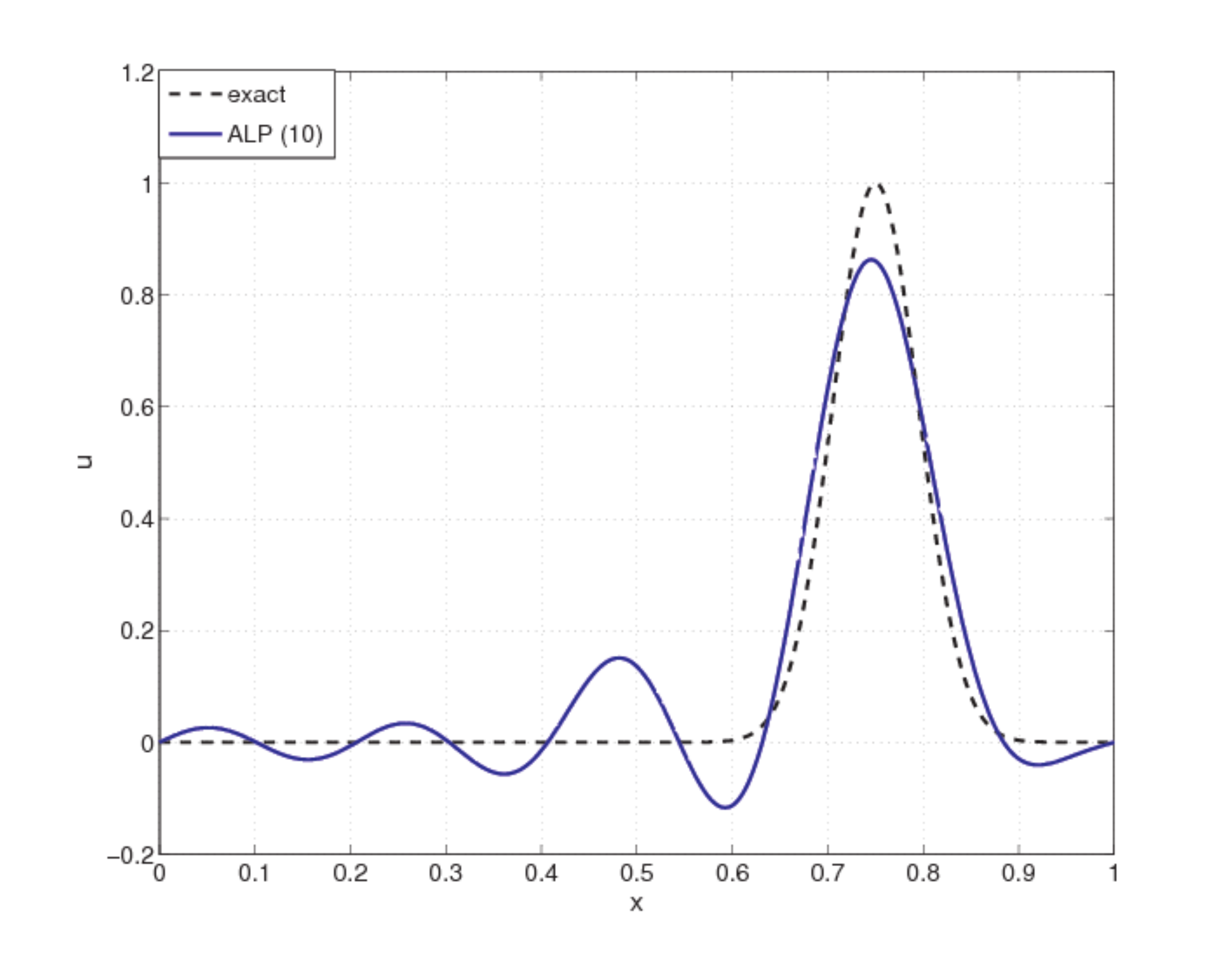} \\[-12pt] (a) \\
\includegraphics[height=7.5cm]{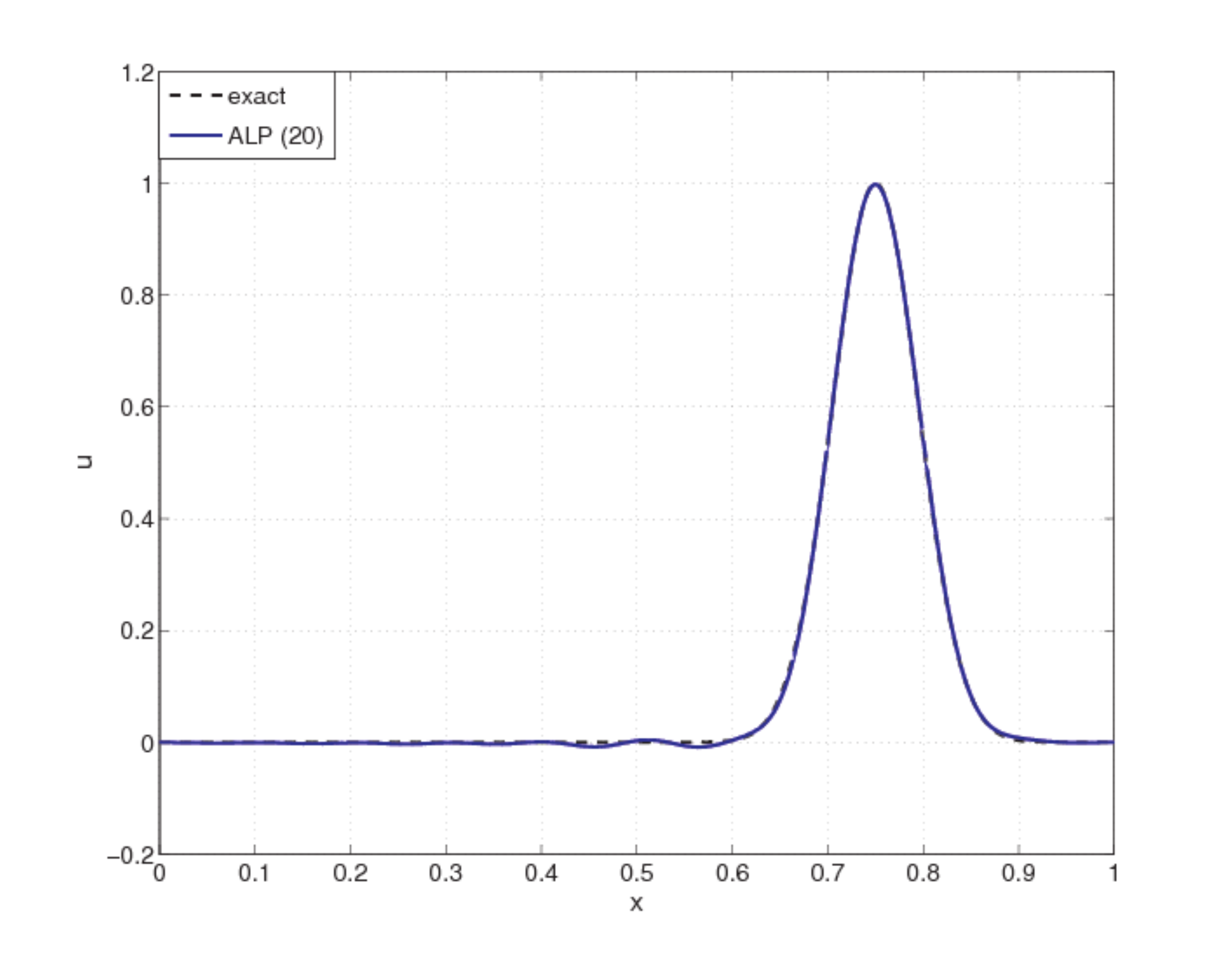}\\\\[-12pt] (b) 
\end{tabular}}}
\caption{Comparison between the exact solution and the reconstruction of the ROM solution for: (a) $N_{M} = 10$, (b) $N_M=20$. The parameter $\chi=150$ is the same for both the simulations, $\delta t = 1/256$. }
\label{fig::LinWaveExactVsALP}
\end{figure}
In Fig.\ref{fig::LinWaveExactVsALP} the comparison between the analytical solution at final time and the reconstruction of the ROM solution in the FEM space is shown. While with $N_M = 10$ modes the solution is not precise and there are large oscillations (Fig.\ref{fig::LinWaveExactVsALP}.(a)), with $N_M=20$ the solution is very accurate (Fig.\ref{fig::LinWaveExactVsALP}.(b)). In both cases the peak position is well captured. 

\begin{rem}
If $M = -cD$, i.e. if $M$ was taken as the discrete form of the known $\mathcal{M}$ operator, then the reduced order form of the PDE would reduce to $\dot{\beta}_i=0$. Thus the reduced order solution would be exactly given by the initial expansion on the modes, advected at a velocity $c$.
\end{rem}

\subsection{Korteweg-de Vries equation}
\label{sec:KdVNumerical}

In this section, our reduced-order method is applied to the Korteweg-de Vries (KdV) equation:
\begin{equation}
\partial_t u + 6 u\partial_x u + \partial^3_x u = 0.
\label{eq:KdVAnalytical}
\end{equation}
The KdV equation is a classical example of integrable system, and a Lax pair is known in closed-form (see section~\ref{sec:kdv-SCSA}). Following the ALP algorithm, this knowledge is not used in the reduced-order model. The expansion of $u$ is injected into the Eq.\eqref{eq:KdVAnalytical} expressed in conservative form, leading to:
\begin{equation}
\sum_{i=1}^{N_M} \dot{\beta}_i \phi_i + \beta_i \partial_t\phi_i + 3 \sum_{i,j=1}^{N_M} \beta_i \beta_j \partial_x(\phi_i \phi_j) + \sum_{i=1}^{N_M}\beta_i \partial^3_x \phi_i = 0. 
\end{equation}
Then the spectral problem is used to simplify this expression. Two strategies may be adopted: either eliminate the quadratic term or transform the third order derivative into a quadratic term. The first strategy would introduce an extra third order tensor, thus increasing the computational cost of the ODE system.
We therefore adopt the second option:
\begin{equation}
-\partial^2_x \phi_i - \chi\sum_{j=1}^{N_M} \beta_j \phi_i\phi_j = \lambda_i \phi_i \ \Rightarrow \chi \sum_{j=1}^{N_M} \beta_j \partial_x (\phi_i\phi_j) = - \lambda_i\partial_x\phi_i - \partial^3_x\phi_i,
\end{equation}
which gives:
\begin{equation}
\sum_{j=1}^{N_M} \dot{\beta}_j \phi_j + \beta_j \partial_t\phi_j - \frac{3}{\chi}\sum_{j=1}^{N_M} (\lambda_j\partial_x\phi_j + \partial^3_x\phi_j)\beta_j + \sum_{j=1}^{N_M}\beta_j\partial^3_x\phi_j = 0.
\end{equation}
This equation is projected onto the eigenfunctions and the evolution of the coefficients $\beta_i$ is obtained:
\begin{equation}
\dot{\beta_i} + \sum_{j=1}^{N_M} M_{ij} \beta_j - \frac{3 }{\chi}\sum_{j=1}^{N_M}\lambda_j D_{ij} \beta_j  + \left(1-\frac{3}{\chi} \right) \sum_{j=1}^{N_M} D^{(3)}_{ij}\beta_j = 0,  
\end{equation}
where $D^{(3)}_{ij}:=\langle \partial^3_x \phi_j, \phi_i\rangle = - \langle \partial^2_x \phi_j,\partial_x\phi_i \rangle$, whose time evolution is governed by $\dot{D}^{(3)} + [D^{(3)},M] = 0$, since it is a linear time independent operator.

The propagation of a one-soliton and of a three-soliton are considered. 
In order to quantify the discrepancy between the analytical solution and the reconstruction, two error indicators are used: the time average and the maximum value over time of $\varepsilon_{L2}(t)$, defined as
\begin{equation}
\varepsilon_{L2}^2 (t) := \frac{\int_{\Omega}(u - u_{ALP})^2\ d\Omega}{\int_{\Omega} u^2\ d\Omega}.
\end{equation}

\paragraph{One-soliton propagation}
\label{sec:kdv1-eigfun}
The exact solution reads:
\begin{equation}
u(x,t) = \frac{\beta}{2} \mathrm{sech}^2\left(\frac{\beta^{1/2}}{2}(x-\beta t - x_0) \right),
\label{KdV_OneSolExact}
\end{equation}
with $\beta = 4$. The final time is set to $T_{max} = 5.0$. 

The modes were extracted by using the initial condition only $u_0=u(x,0)$, setting $\chi = 1$. The Schr\"{o}dinger spectral problem was discretized in a space of $N_h=500$ piecewise linear functions.
\begin{table}
\begin{center}
\begin{tabular}{cccc}
\hline
& $N_M$ & $\overline{\varepsilon_{L2}}$ & $\max_t(\varepsilon_{L2})$ \\
\hline
& $26$ & $0.2176$ & $0.3272$ \\
& $28$ & $0.2007$ & $0.3168$ \\
& $30$ & $0.1388$ & $0.1869$ \\
& $32$ & $0.0916$ & $0.1394$ \\
& $34$ & $0.0486$ & $0.0637$ \\
& $36$ & $0.0370$ & $0.0578$ \\
\hline
\end{tabular}\caption{Error indicators for the KdV one-soliton test case (Section \ref{sec:kdv1-eigfun}) as a function of the number of modes used to discretize the equations: first column $N_M$ is the number of modes, the second and the third one the time average and the maximum of the $L^2$ error. To be compared to Table \ref{table::KdV1-solit}.}\label{table::KdV1-eigfun}
\end{center}
\end{table}
In Table \ref{table::KdV1-eigfun}, the error indicators are shown as function of the number of modes used to discretize the system. 
\begin{figure}
\centerline{\hbox{\begin{tabular}{cc}
\includegraphics[height=7.5cm]{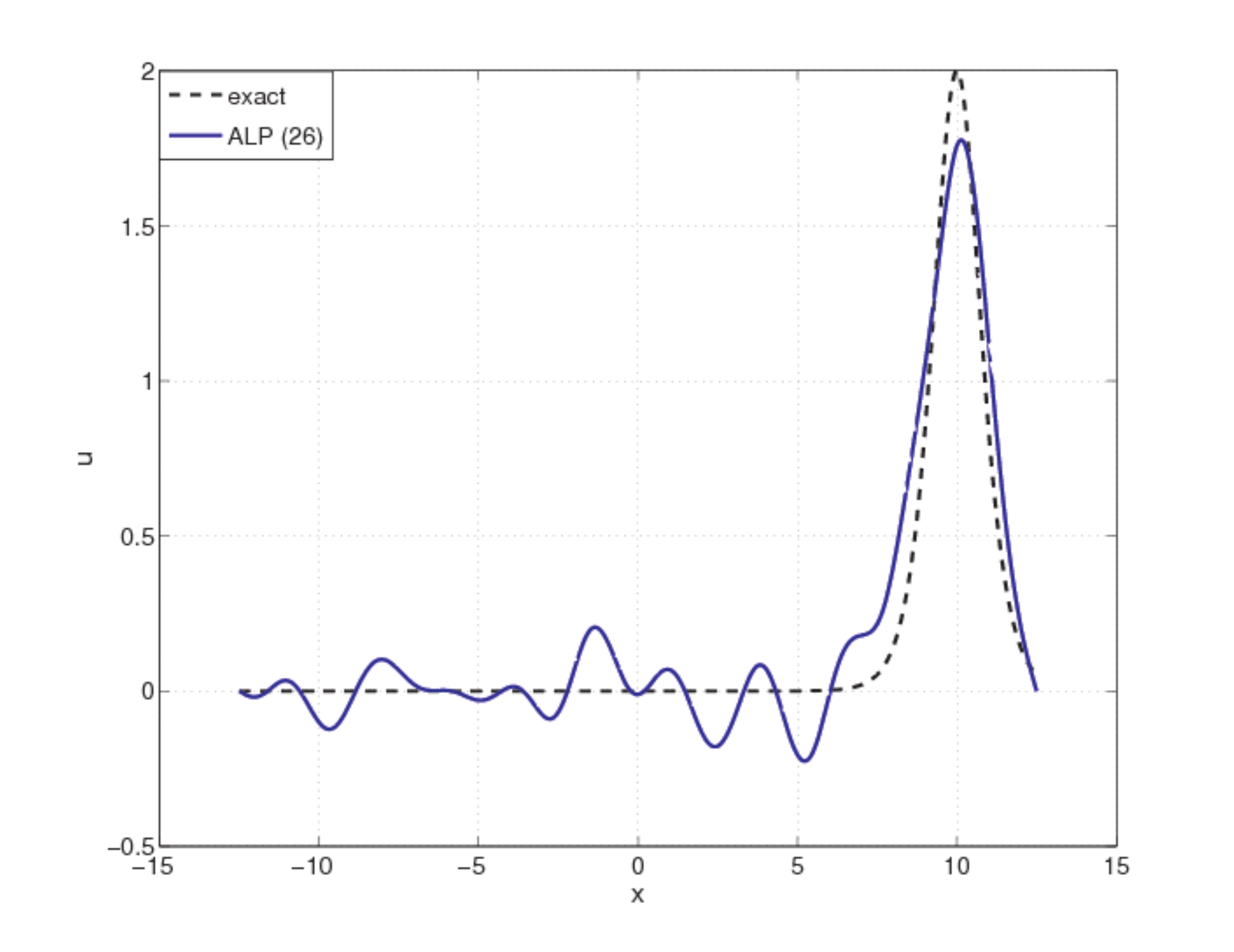} \\ (a) \\
\includegraphics[height=7.5cm]{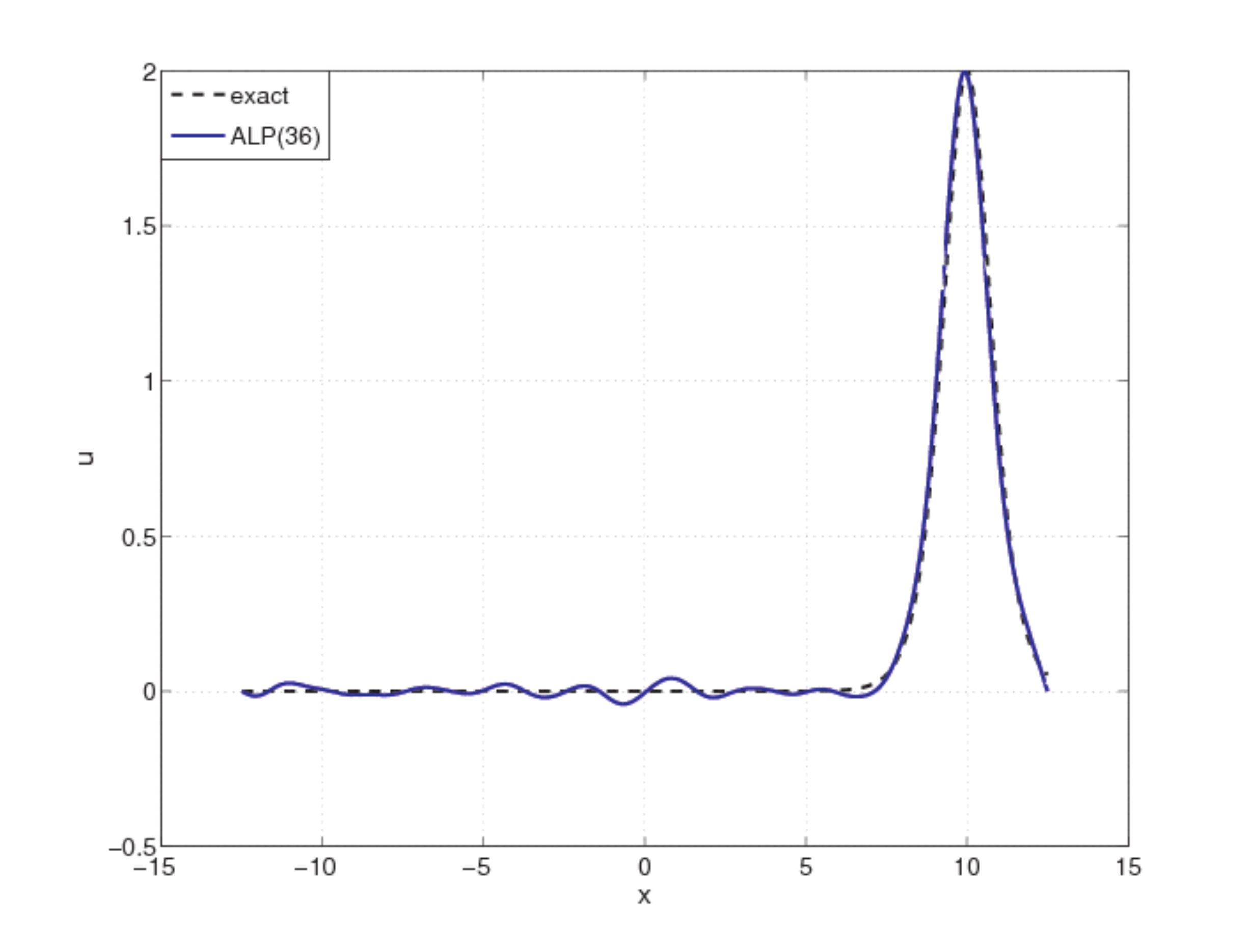}\\ (b) \\
\end{tabular}}}
\caption{Comparison at final time between the exact solution and the solution obtained by integrating ALP with (a) $N_M=26$, (b) $N_M=36$.}
\label{fig::KdV1-eigfun}
\end{figure}
In Fig.\ref{fig::KdV1-eigfun} a comparison between the reconstruction of the ROM and the analytical solution is proposed, for $N_M=26$ and $N_M=36$, $t=T_{max}$. When using $N_M=26$ some errors in the shape and in the amplitude are still present, while, when $N_M=36$, the profile motion is well captured, the error being only concentrated is small oscillations behind the wave.
The main difficulty in integrating this test case is due to the large distance travelled by the wave, which is characterized by a relatively sharp profile. 
However, it is worth noting that the error is mainly due to the reconstruction (post-processing) stage. As said in Section~\ref{sec:reducedToFull}, a better reconstruction scheme will be the object of future works.

\begin{figure}
\centerline{\hbox{\begin{tabular}{cc}
\includegraphics[height=7cm]{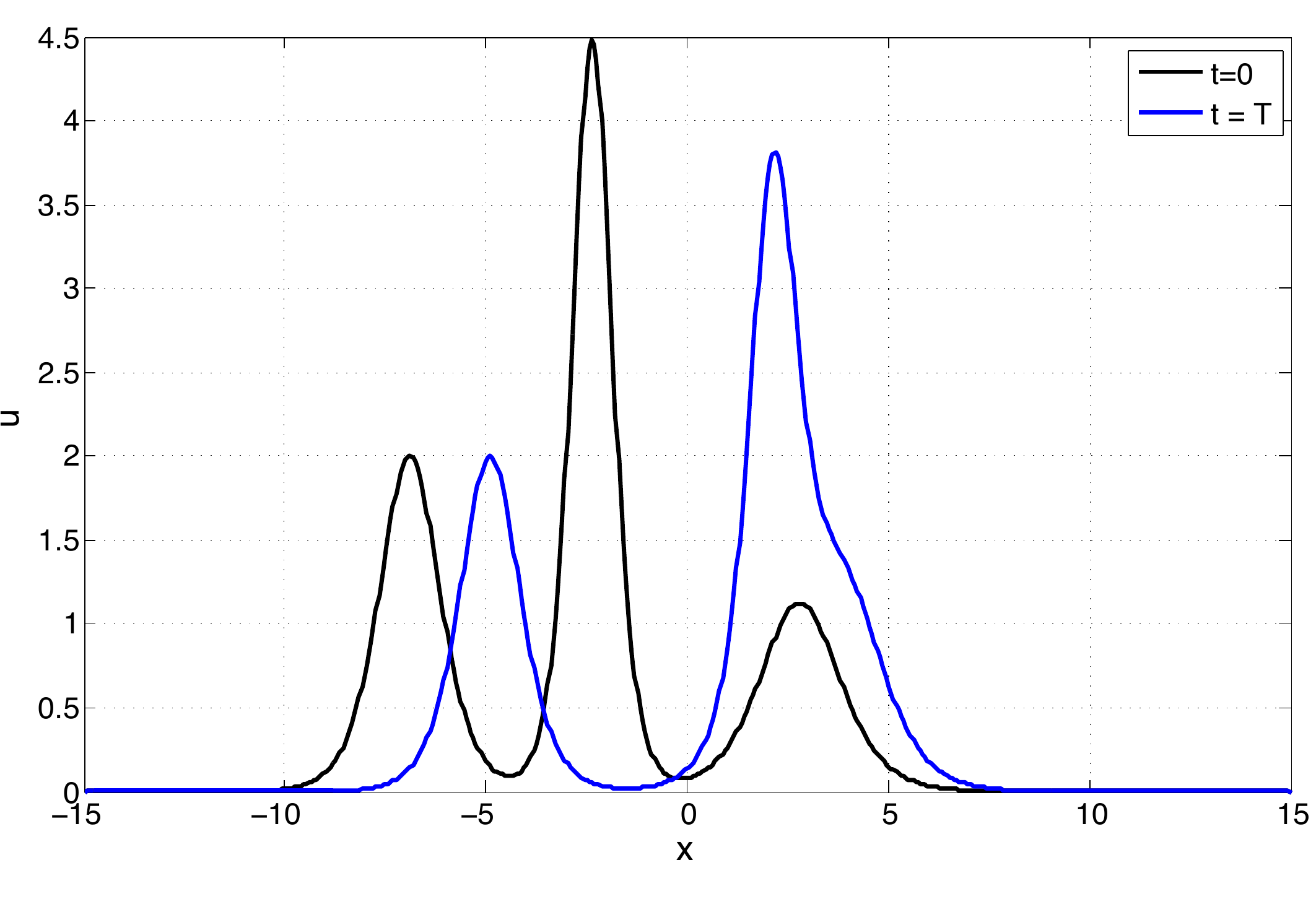} \\ (a) \\
\includegraphics[height=7.5cm]{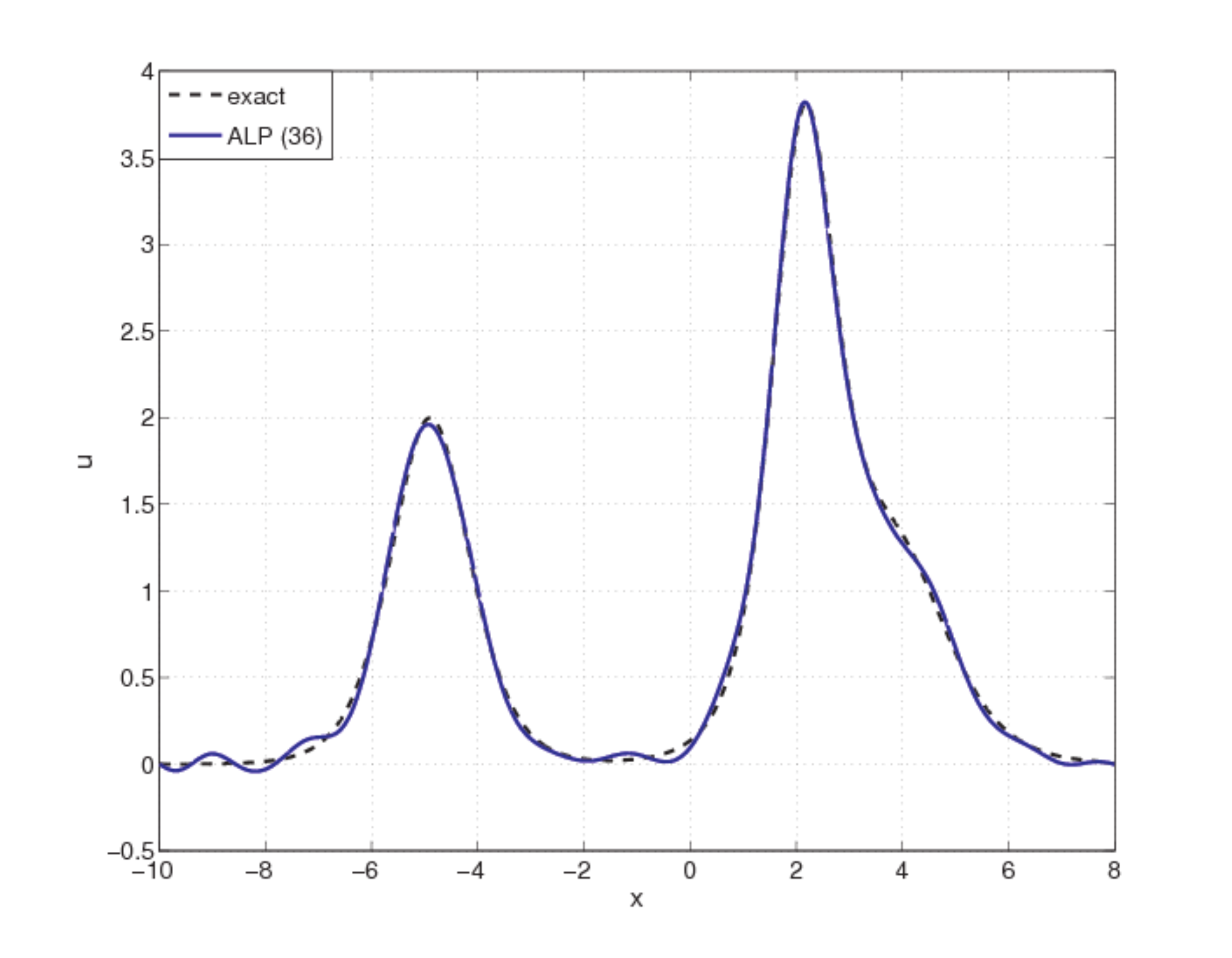}\\ (b)
\end{tabular}}}
\caption{(a) Initial and final configuration for the three-soliton solution (b) Comparison at $t=T_{max}$ between the exact solution and the solution obtained by integrating ALP with $N_M=36$.}
\label{fig::KdV_3-eigfun}
\end{figure}

\paragraph{Three-soliton propagation}
A three-soliton propagation is taken as an example of  $n$ solitary interacting waves. The reference solution, shown in Fig.\ref{fig::KdV_3-eigfun}.(a), has been generated by considering the Gelfand-Levitan-Marchenko equation, that, when solved for the KdV equation provides:
\begin{equation}
u(x,t) = -2\partial^2_{x} \log (\det (I+A(x,t))),
\label{threeSolExact}
\end{equation} 
where $A\in\mathbb{R}^{n\times n}$ is the interacting matrix, written in terms of the scattering data \cite{Ablowitz_1}. In particular:
\begin{equation}
A_{mn}(x,t) = \frac{c_m c_n}{k_m + k_n} \exp\left\lbrace (k_m+k_n)x - (k_m^3+k_n^3)t \right\rbrace,
\end{equation}
where $k_m,c_n$ are $2n$ scalar parameters that may be linked to position and speed of solitons (see \cite{Drazin_1}).
For the present case: $c = [5.0\ 10^{-2},1.5\ 10^{-1},1.0\ 10^1]$, $k = [1.0, 1.5, 1.75]$, $x\in(-15,15)$ and $t\in(0,0.5)$.
This setting is challenging because of the interaction of the waves: at final time, two of them are fused together (Fig.\ref{fig::KdV_3-eigfun}).

The error indicators are computed as function of the number of modes. The results are shown in Table \ref{table::KdV3-eigfun}.
\begin{table}
\begin{center}
\begin{tabular}{cccc}
\hline
& $N_M$ & $\overline{\varepsilon_{L2}}$ & $\max_t(\varepsilon_{L2})$ \\
\hline
& $28$ & $0.0896$ & $0.0991$ \\
& $32$ & $0.0558$ & $0.0690$ \\
& $36$ & $0.0249$ & $0.0295$ \\
& $40$ & $0.0184$ & $0.0215$ \\
& $44$ & $0.0138$ & $0.0180$ \\
& $48$ & $0.0081$ & $0.0121$ \\
\hline
\end{tabular}
\caption{Error indicators for the KdV three-soliton test case as a function of the number of modes used to discretize the equations: first column $N_M$ is the number of modes, the second and the third one the time average and the maximum of the $L^2$ error.}\label{table::KdV3-eigfun}
\end{center}
\end{table}
The qualitative behavior of the scheme is good at reproducing the dynamics of the waves interaction, as it can be seen in Fig.\ref{fig::KdV_3-eigfun}.(b). The error is mainly due to small oscillations arising in the flat part of the domain.
The peak positions of the waves, as well as their amplitude is correct.

\subsection{Fisher-Kolmogorov-Petrovski-Piskunov equation}
In this section the Fisher-Kolmogorov-Petrovski-Piskunov (FKPP) equation is considered as an example of non-isospectral flow equation, in a finite domain, with homogeneous Dirichlet or Neumann boundary conditions. This is a case for which \eqref{eq:lax-isospec} is \emph{a priori} not satisfied. 


\subsubsection{1D FKPP with homogeneous Dirichlet boundary conditions}
\label{sec:fkpp-1d}

The equation reads:
\begin{equation}
	\left\{
\begin{split}
& \partial_t u = \partial_x^2 u + \nu u(1-u), \ \mbox{ in }\ \Omega = [0,1],\\
& u = 0, \ \mbox{ on }\ \partial \Omega.
\end{split}
\right . 
\label{FKPP_Eq1D}
\end{equation}
For the present case $\nu = 10^3$, the space domain is $[0,1]$ and $N_h = 250$. The time domain is $[0,7.5\cdot 10^{-3}]$ and $100$ integration points are taken. The reference solution is obtained by discretizing in space by means of piecewise linear functions and by using a mixed implicit-explicit scheme in time: the linear diffusion part of the equation is discretized by a Cranck-Nicolson scheme, the nonlinear term by an explicit second order Adams-Bashforth scheme, with $\delta t = 7.5 \ 10^{-5}$.
\begin{figure}
\centerline{\hbox{\begin{tabular}{cc}
\includegraphics[height=7.5cm]{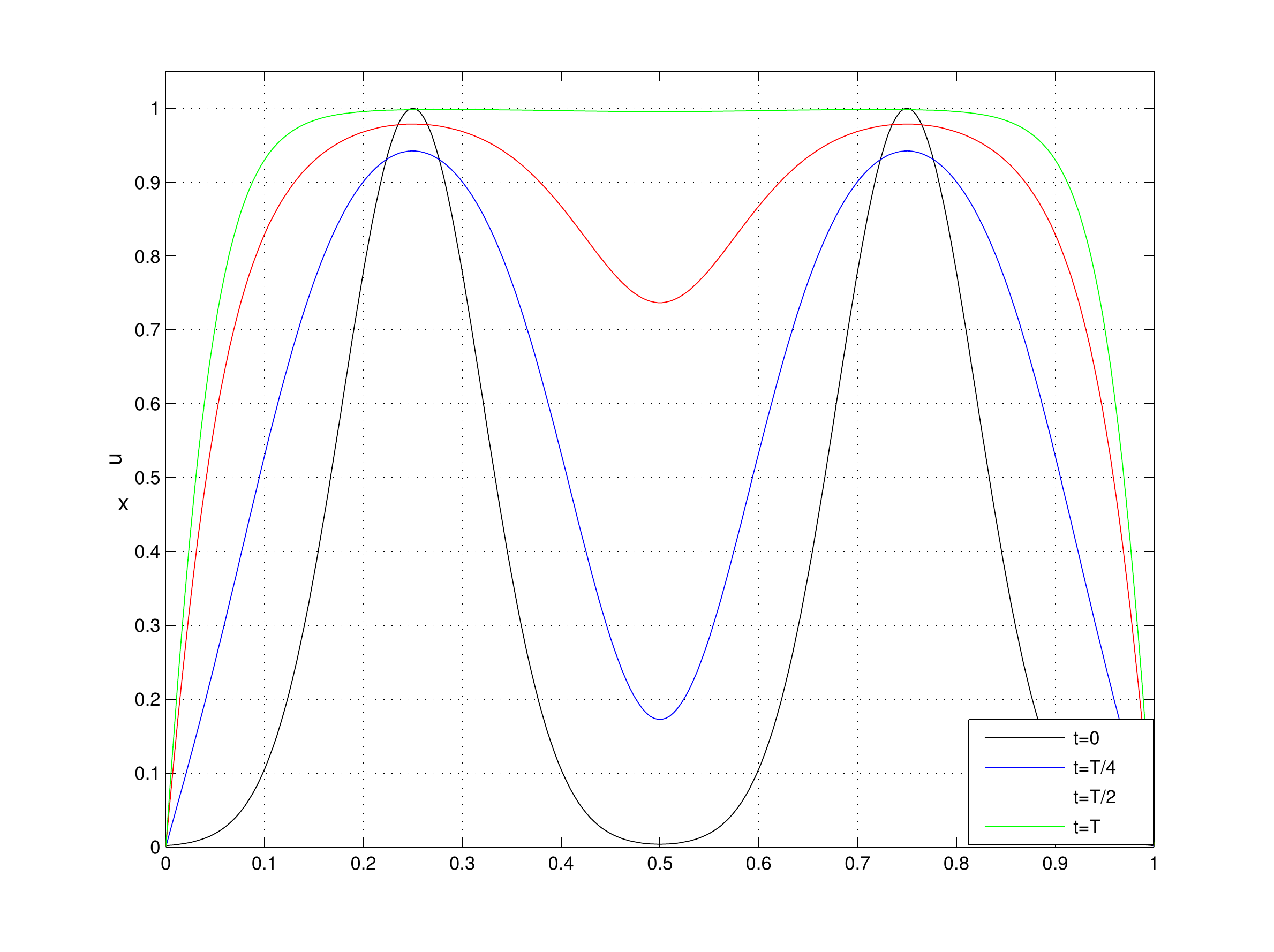} \\ (a) \\
\includegraphics[height=7cm]{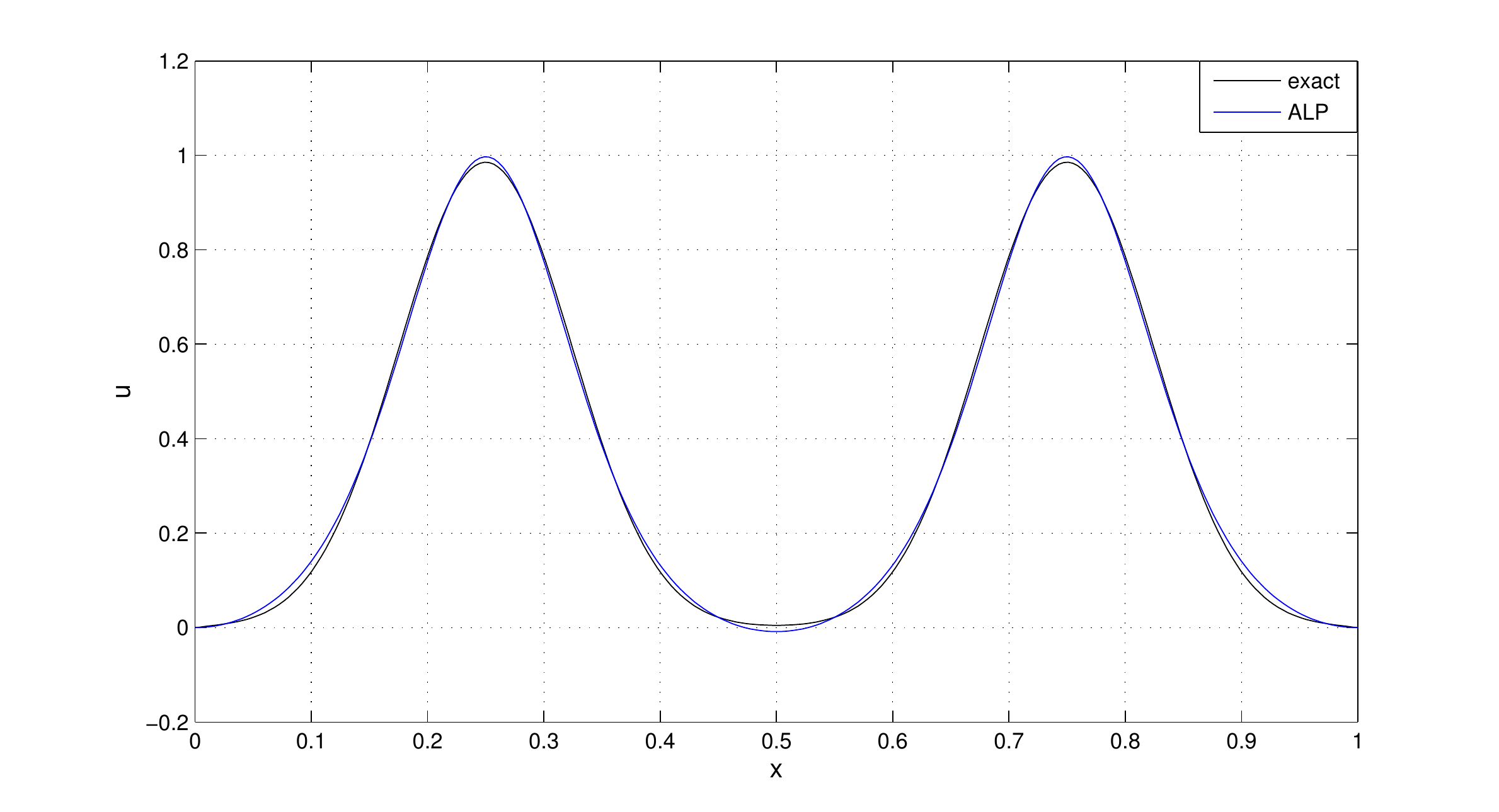}\\(b)
\end{tabular}}}
\caption{(a) Plot of the reference solution of \eqref{FKPP_Eq1D}, at different times; (b) Comparison between the exact solution and the reconstructed one at initial time by using four modes.}
\label{FKPP_Sol}
\end{figure}
The initial solution is given by $u_0 = \exp\left(-10^2(x-0.25)^2\right)+\exp\left(-10^2(x-0.75)^2\right)$. 

Inserting the modal approximation of the solution in the equation and using $\Delta \phi_i = -(\lambda_i + \chi u)\phi_i$, the following holds:
\begin{equation}
	\label{eq:beta-fkpp}
	\partial_t \beta_i + \sum_{j=1}^{N_M} M_{ij}\beta_j =  (\nu-\lambda_i) \beta_i - (\chi + \nu) \sum_{j,k=1}^{N_M} \mathcal{T}_{ijk}\beta_j \beta_k.
\end{equation}

Identifying \eqref{eq:beta-fkpp} with $\eqref{eq:rom-dynamics}_1$ gives $\eqref{eq:rom-dynamics}_5$, i.e. the relation between $\beta$ and $\gamma$ specific to the FKPP equation:
\begin{equation}
	\label{eq:gamma-fkpp}
	\gamma_i = (\nu-\lambda_i) \beta_i - (\chi + \nu) \sum_{j,k=1}^{N_M} \mathcal{T}_{ijk}\beta_j\beta_k.
\end{equation}
The system to be solved in the reduced space is therefore $\eqref{eq:rom-dynamics}_{1-4}$, \eqref{eq:gamma-fkpp}.

The error indicators considered to investigate the behavior of the ROM are:
\begin{eqnarray}
\overline{\varepsilon}_{L2}^2 := \frac{1}{T_{max}}\int_{0}^{T_{max}} \varepsilon_{L2}^2(t) \ dt =  \frac{1}{T_{max}}\int_0^{T_{max}} \frac{\int_{\Omega}(u - u_{ALP})^2\ d\Omega}{\int_{\Omega} u^2\ d\Omega} \ dt, \\
\varepsilon_{T_{max}}^2  := \varepsilon_{L2}^2(t=T_{max}).
\end{eqnarray}
In Table \ref{table::FKPP1DErr}, the values of the error indicators are written as a function of the number of modes used. The performance of the method is overall satisfactory. 
\begin{table}
\begin{center}
\begin{tabular}{cccc}
\hline
& $N_M$ & $\overline{\varepsilon}_{L2}$ & $\varepsilon_{T_{max}}$ \\
\hline
& $6$ & $0.1722$ & $0.2218$ \\
& $8$ & $0.0522$ & $0.0747$ \\
& $10$ & $0.0304$ & $0.0458$ \\
& $12$ & $0.0163$ & $0.0279$ \\
& $14$ & $0.0097$ & $0.0168$ \\
& $16$ & $0.0059$ & $0.0105$ \\
\hline
\end{tabular}
\caption{Error indicators for the 1D FKPP test case as a function of the number of modes used to discretize the equations: first column $N_M$ is the number of modes, the second and the third one the average error in $L^2$ norm and the error at final time.}\label{table::FKPP1DErr}
\end{center}
\end{table}
As done for the linear advection equation, the Frobenius norm indicator~\eqref{eq::frobErr} is monitored. In Fig.\ref{fig::Frobenius}.(b) the time average and the maximum of this error indicator as a function of time are shown.

\subsubsection{2D FKPP with homogeneous Neumann boundary conditions}

The bidimensional FKPP equation reads:
\begin{equation}
	\left\{
\begin{split}
\partial_t u = \Delta u + \nu u (1-u), \mbox{ in } \Omega, \\
\partial_n u = 0, \mbox{ on } \partial\Omega,
\end{split}
\right .
\label{FKPP_Eq2D}
\end{equation}
where $\Omega$ is a bounded domain of  $\mathbb{R}^2$. 

\paragraph{Unit square geometry}
In this test case, $\Omega$ is a unit square. The number of degrees of freedom is about $N_h = 5700$. The logistic coefficient is $\nu=50$, the final time $T_{max}=5\ 10^{-2}$ and $\delta t = 5\ 10^{-4}$, so that $100$ time iterations are performed. The same time step was considered for the ROM integration.  

\begin{figure}
\centerline{\hbox{\begin{tabular}{cc}
\includegraphics[height=6.0cm,width=7.0cm]{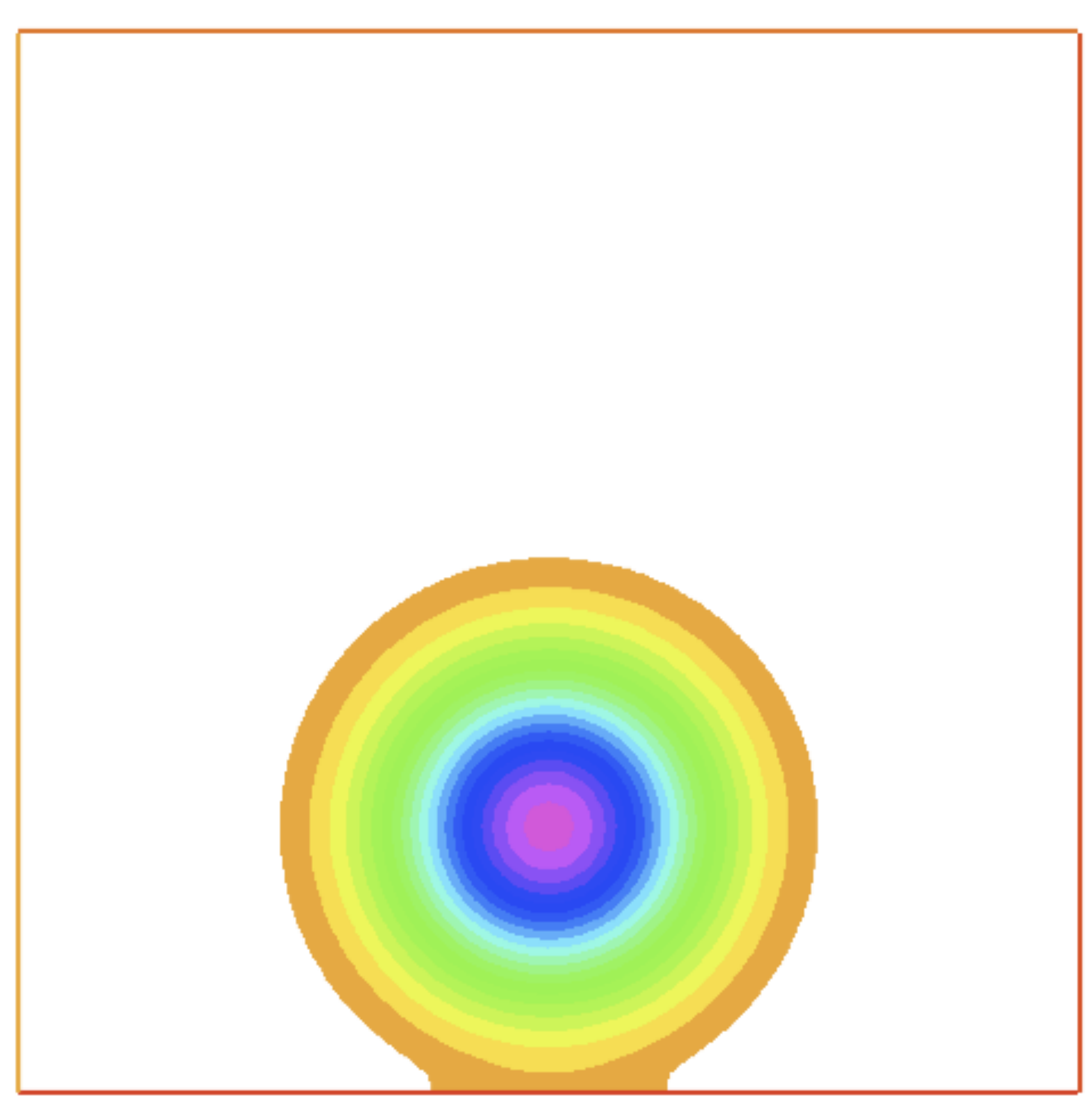} &
\includegraphics[height=6.0cm,width=7.0cm]{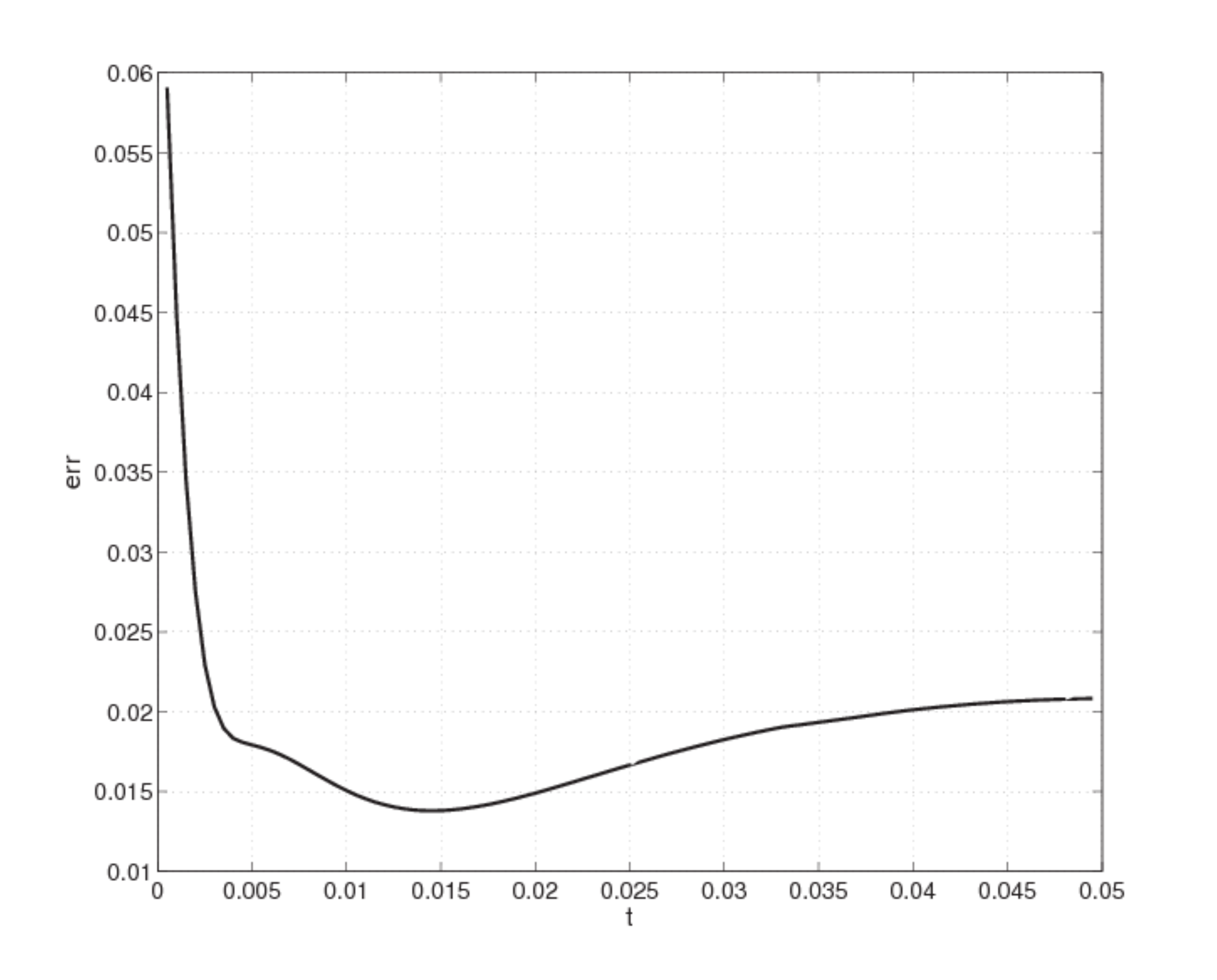}\\
\hspace{0.2cm} (a) & \hspace{0.2cm}(b) 
\end{tabular}}}
\caption{2D FKPP equation on a square: (a) Initial datum, (b) $L^2$ error as a function of time when considering $N_M=40$ modes and $\chi=25$.}
\label{fig::FKPP2D_u0}
\end{figure}

\begin{figure}
\centerline{\hbox{\begin{tabular}{ccc}
\includegraphics[height=5.0cm,width=5.5cm]{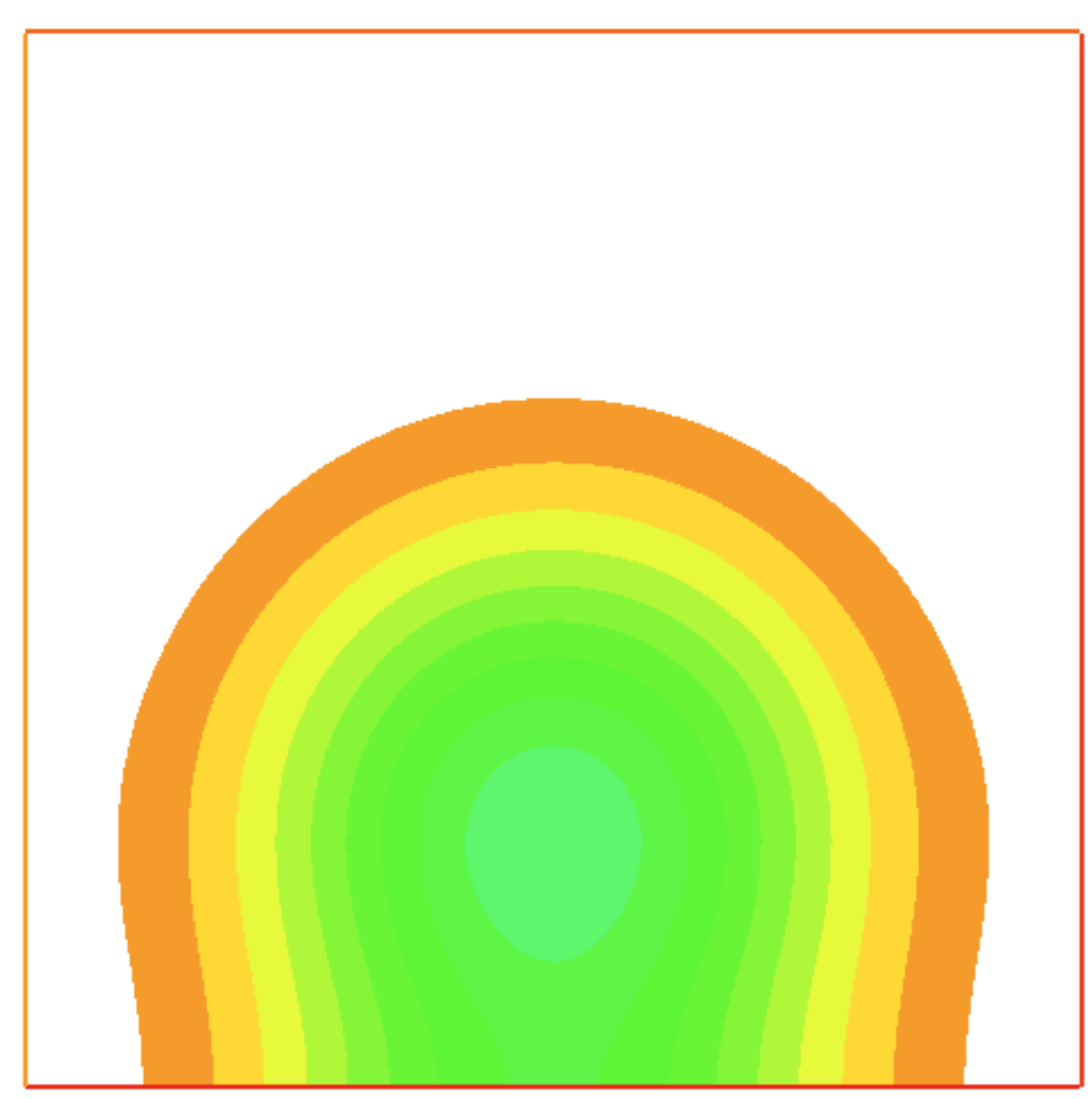} &
\includegraphics[height=5.0cm,width=5.5cm]{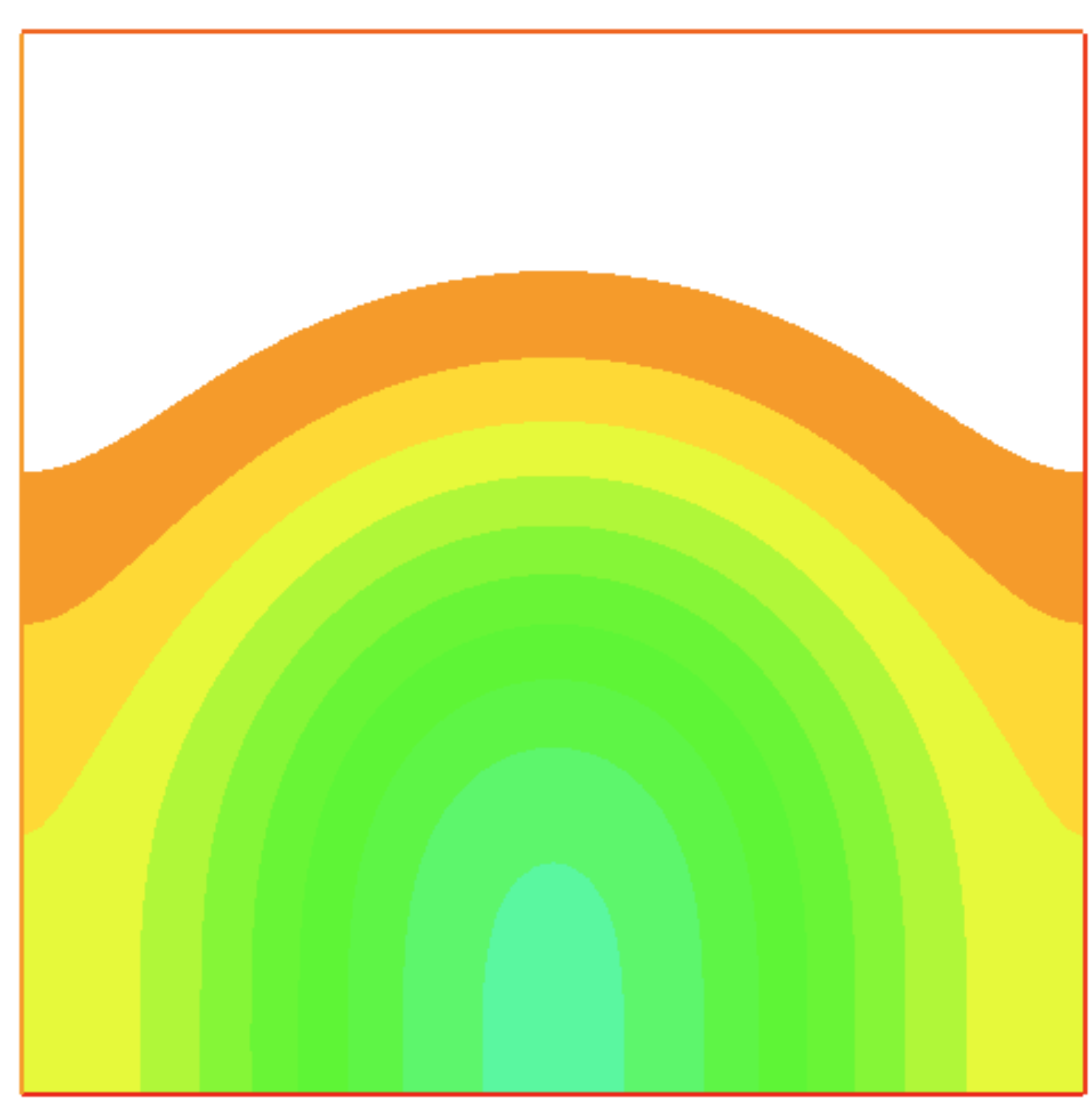} &
\includegraphics[height=5.0cm,width=5.5cm]{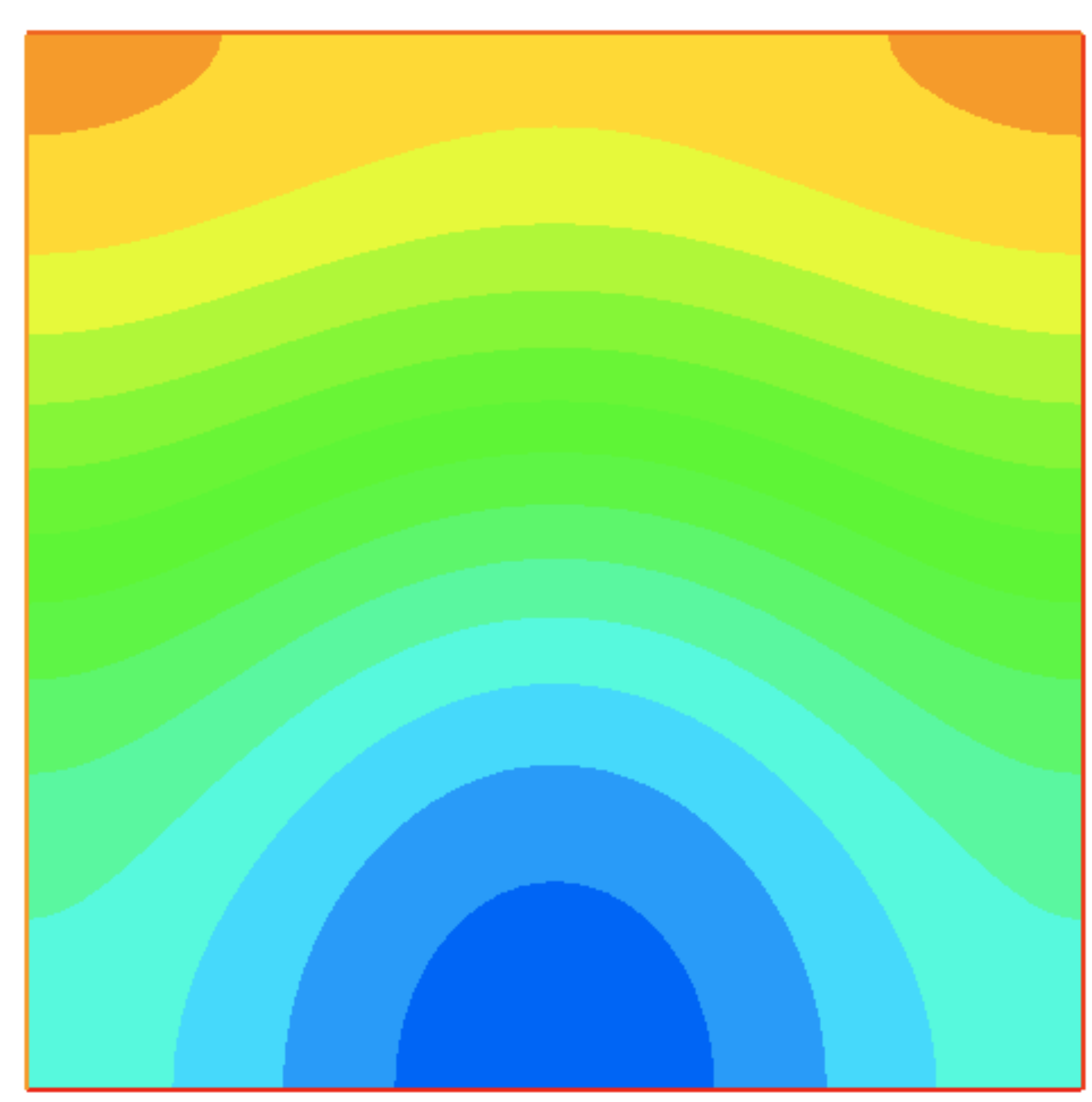}\\
\hspace{0.1cm} (a) & \hspace{0.1cm}(b) & \hspace{0.1cm}(c)
\end{tabular}}}
\caption{2D FKPP reference solution at different times: (a) $t=T_{max}/4$, (b) $t=T_{max}/2$, (c) $t=T_{max}$.}
\label{fig::FKPP2D_ref}
\end{figure}

The initial datum is $u_0(x,y) = \exp\left(-50 ((x-0.5)^2 + (y-0.25)^2) \right)$, whose isovalues are represented in Fig.\ref{fig::FKPP2D_u0}.(a). 
The solution $u$ gets closer to the lower boundary (see for instance Fig.\ref{fig::FKPP2D_ref}.(a) ) in an initial phase, then a front tends to form and propagates upwards (as it is represented in Fig.\ref{fig::FKPP2D_ref}.(b-c)). 
In Fig.\ref{fig::FKPP2D_u0}.(b) the $L^2$ error of the reconstruction with $N_M=40$ modes is shown as a function of time. The symmetry is not perfectly respected at the discrete level because the FEM mesh is unstructured and not symmetric

\begin{figure}
\centerline{\hbox{\begin{tabular}{ccc}
\includegraphics[height=5.0cm,width=5.5cm]{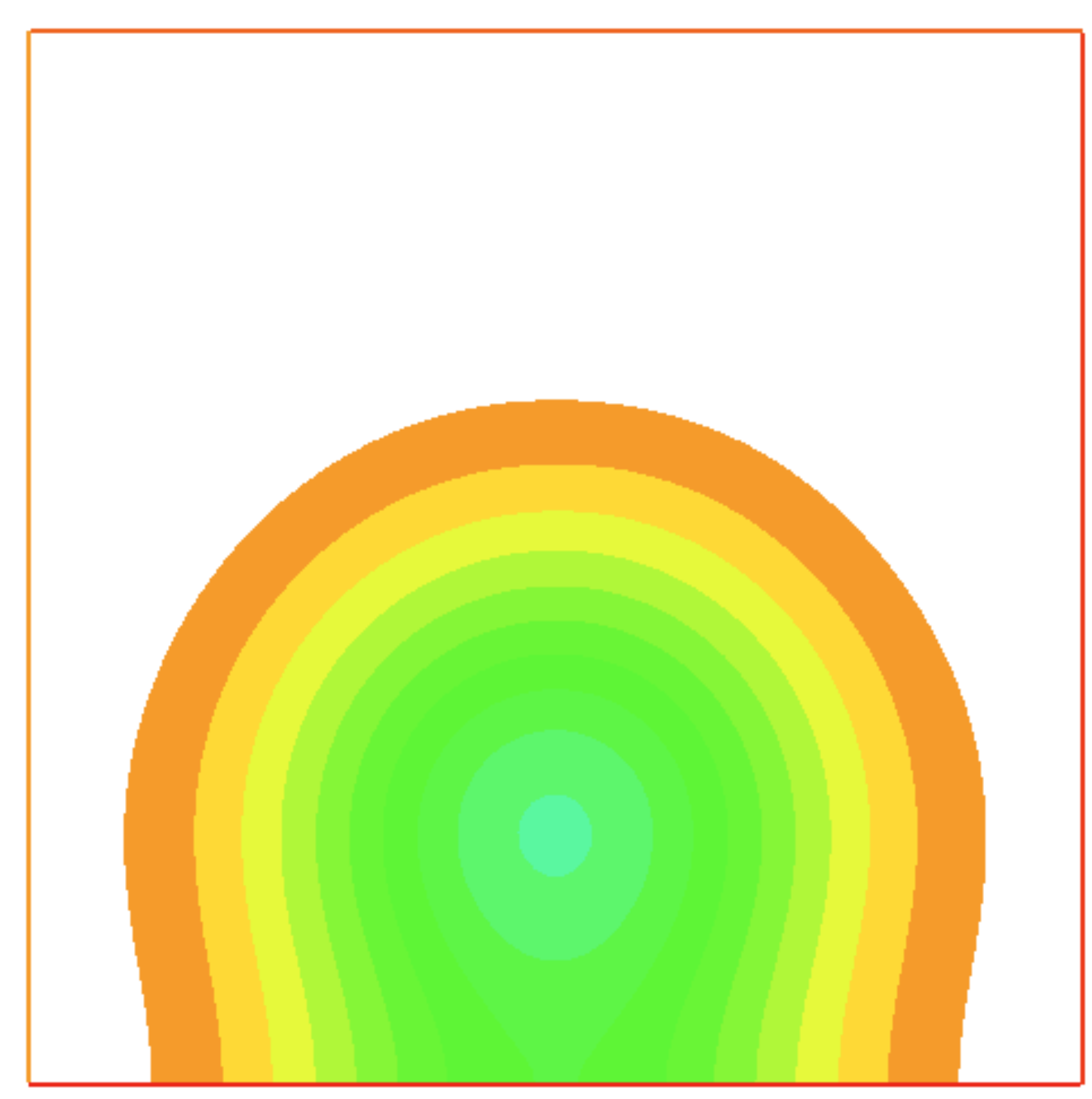} &
\includegraphics[height=5.0cm,width=5.5cm]{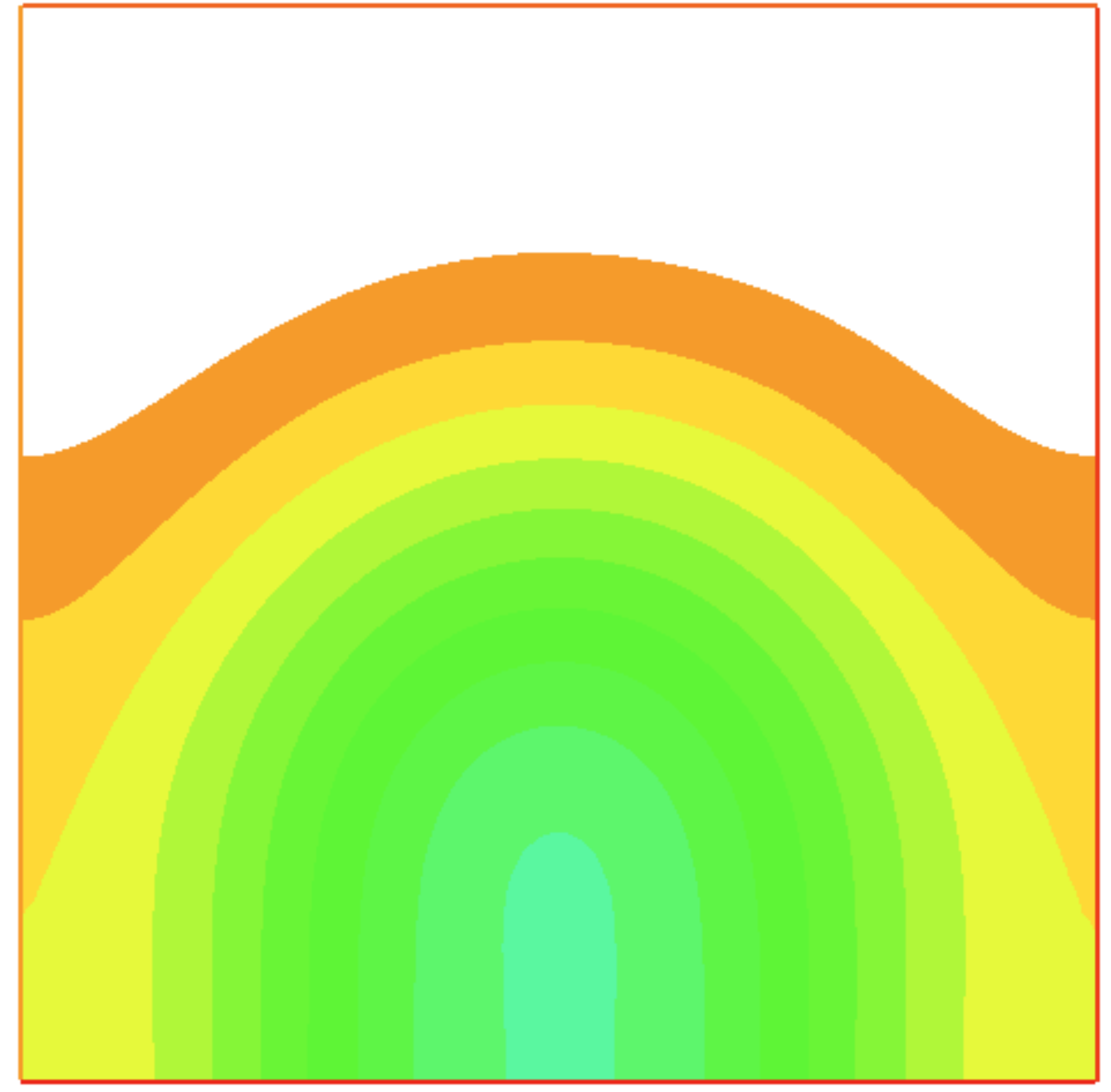} &
\includegraphics[height=5.0cm,width=5.5cm]{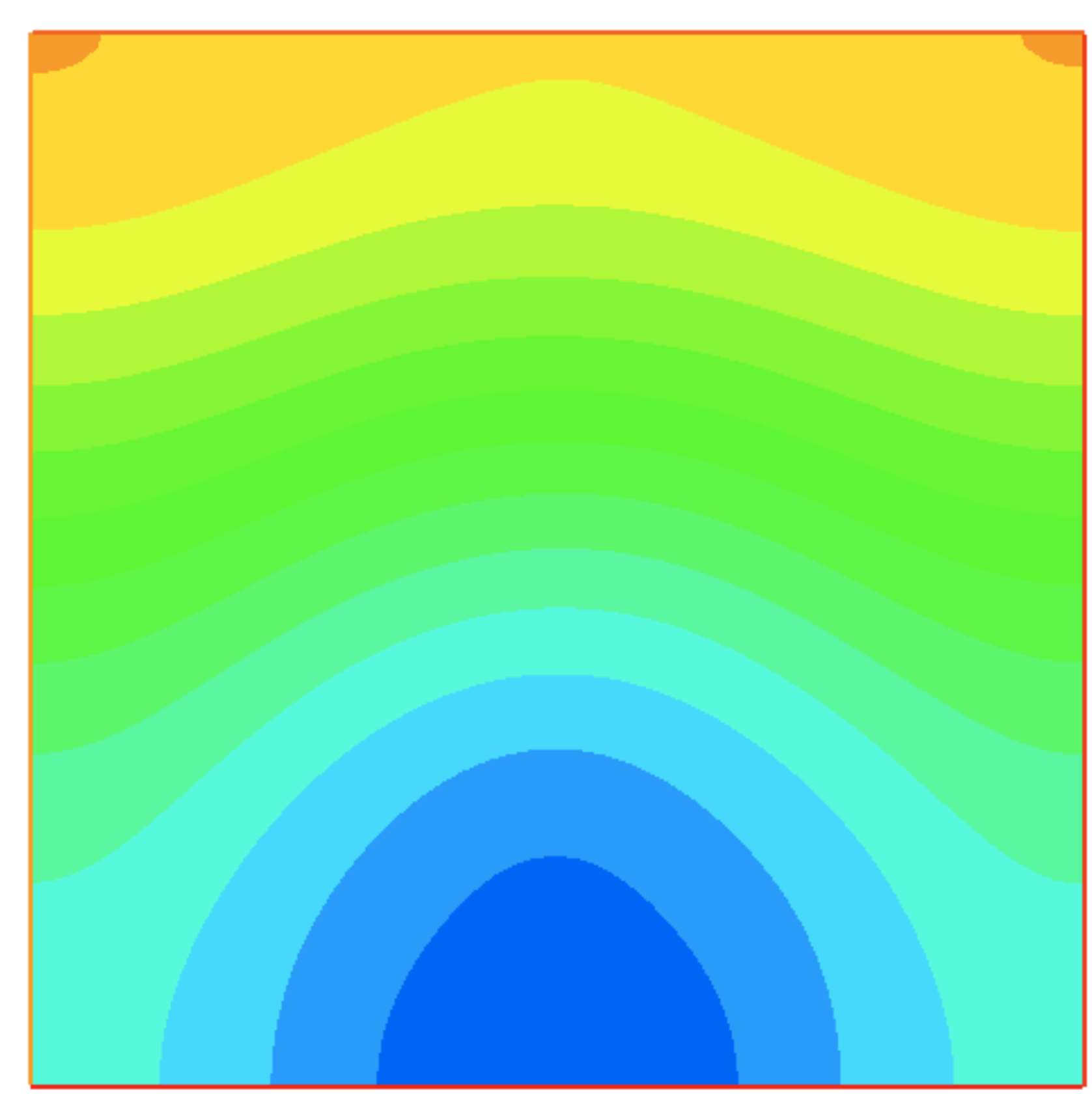}\\
\hspace{0.1cm} (a) & \hspace{0.1cm}(b) & \hspace{0.1cm}(c)
\end{tabular}}}
\caption{2D FKPP ROM solution, obtained with $N_M=40,\ \chi=25$ at different times: (a) $t=T_{max}/4$, (b) $t=T_{max}/2$, (c) $t=T_{max}$.}
\label{fig::FKPP2D_alp}
\end{figure}

The qualitative behavior of the reconstruction may be judged by comparing Fig.\ref{fig::FKPP2D_ref} (the reference solution) and Fig.\ref{fig::FKPP2D_alp}, which shows the reconstruction, after post-processing, of the solution obtained by ALP when $N_M=40$ and $\chi=25$. The symmetry of the solution is not perfectly recovered, but the dynamical behavior is satisfactory, and the error at final time is reasonable, given that the number of degrees of freedom has been divided by about 150 with respect to the FEM solution. The same error indicators introduced for the 1D case are monitored and the results of the numerical simulations are written in Table \ref{table::FKPP2DErr}.
\begin{table}
\begin{center}
\begin{tabular}{cccc}
\hline
& $N_M$ & $\overline{\varepsilon}_{L2}$ & $\varepsilon_{T_{max}}$ \\
\hline
& $5$ & $0.2152$ & $0.0908$ \\
& $10$ & $0.1059$ & $0.0432$ \\
& $15$ & $0.0837$ & $0.0354$ \\
& $20$ & $0.0432$ & $0.0270$ \\
& $25$ & $0.0241$ & $0.0236$ \\
& $30$ & $0.0203$ & $0.0234$ \\
\hline
\end{tabular}
\caption{Error indicators for the 2D FKPP test case on the unit square as a function of the number of modes used to discretize the equations: first column $N_M$ is the number of modes, the second and the third one the average error in $L^2$ norm and the error at final time.}\label{table::FKPP2DErr}
\end{center}
\end{table}

\begin{figure}
\centerline{\hbox{\begin{tabular}{c}
\includegraphics[height=6.0cm,width=7.5cm]{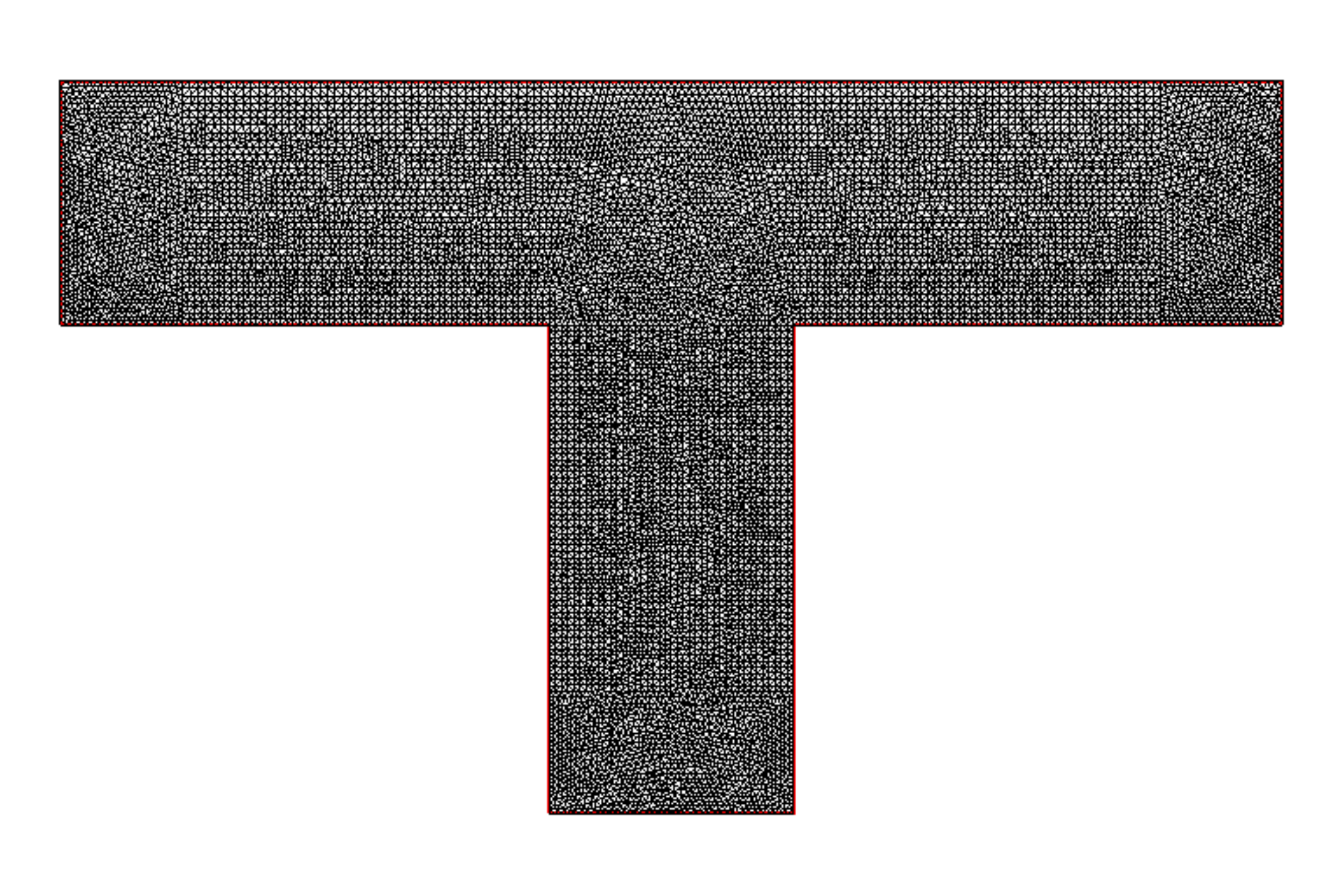}
\end{tabular}}}
\caption{2D FKPP T-shape test case: mesh, with $N_h\approx 11300$ vertices.}
\label{fig::FKPP2D_mesh}
\end{figure}

\paragraph{T-shape geometry}
The same method has been applied to a T-shape geometry in which the front propagates and split (Fig.\ref{fig::FKPP2D_mesh}). 
The direct simulation was performed by using $P1$ finite elements for space discretization: the number of degrees of freedom was $N_h \approx 11300$.
The final time is $T_{max}=0.1364$ and $\delta t = 1.1\ 10^{-3}$, the logistic parameter is set to $\nu=75$ . For the reduced order model, the scattering constant was set to $\chi=1$, the number of modes retained was varied to study the discretization properties, the time step was kept equal to that of the direct simulation.
\begin{figure}
\centerline{\hbox{\begin{tabular}{ccc}
\includegraphics[height=5.0cm,width=5.5cm]{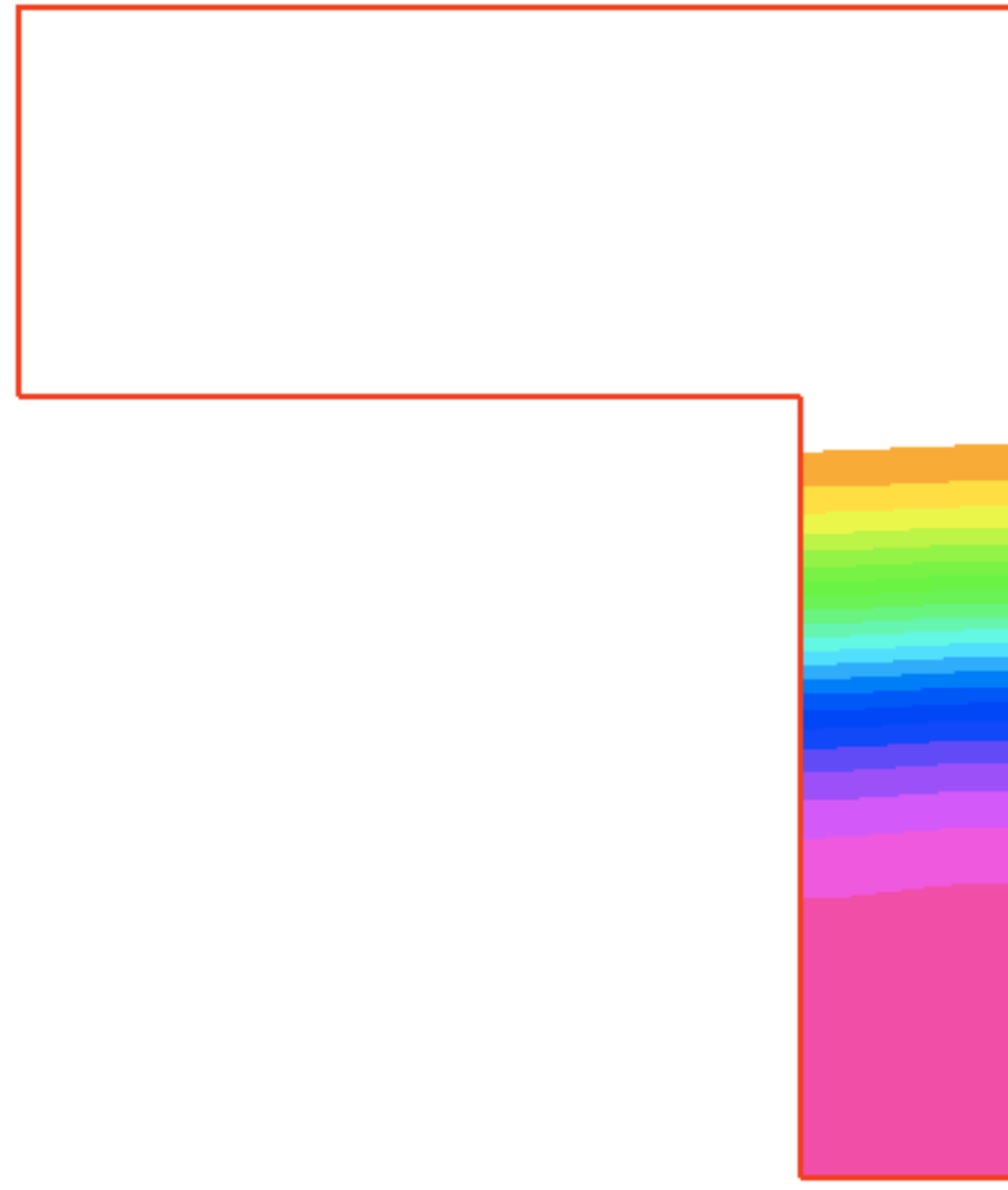} &
\includegraphics[height=5.0cm,width=5.5cm]{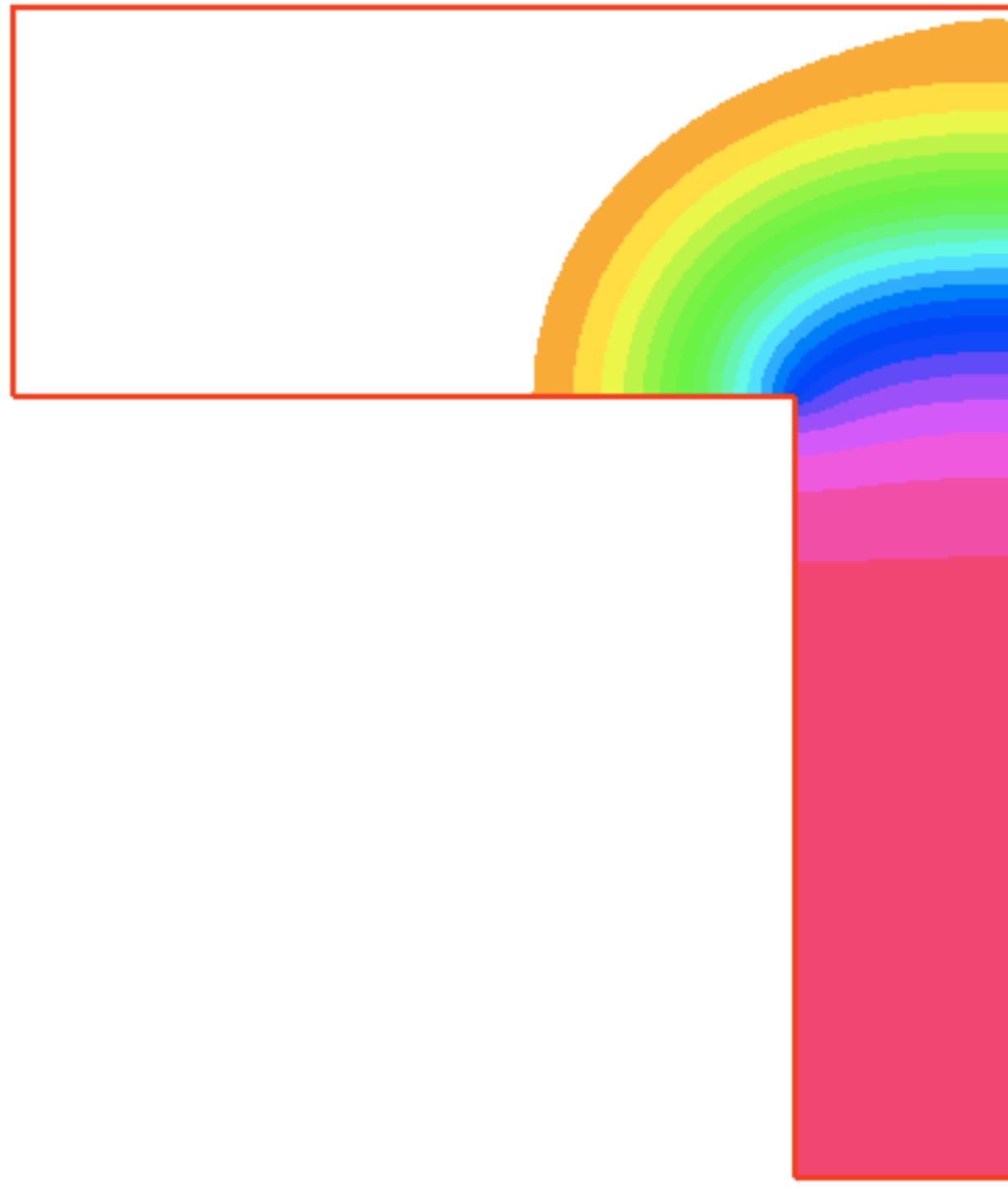} &
\includegraphics[height=5.0cm,width=5.5cm]{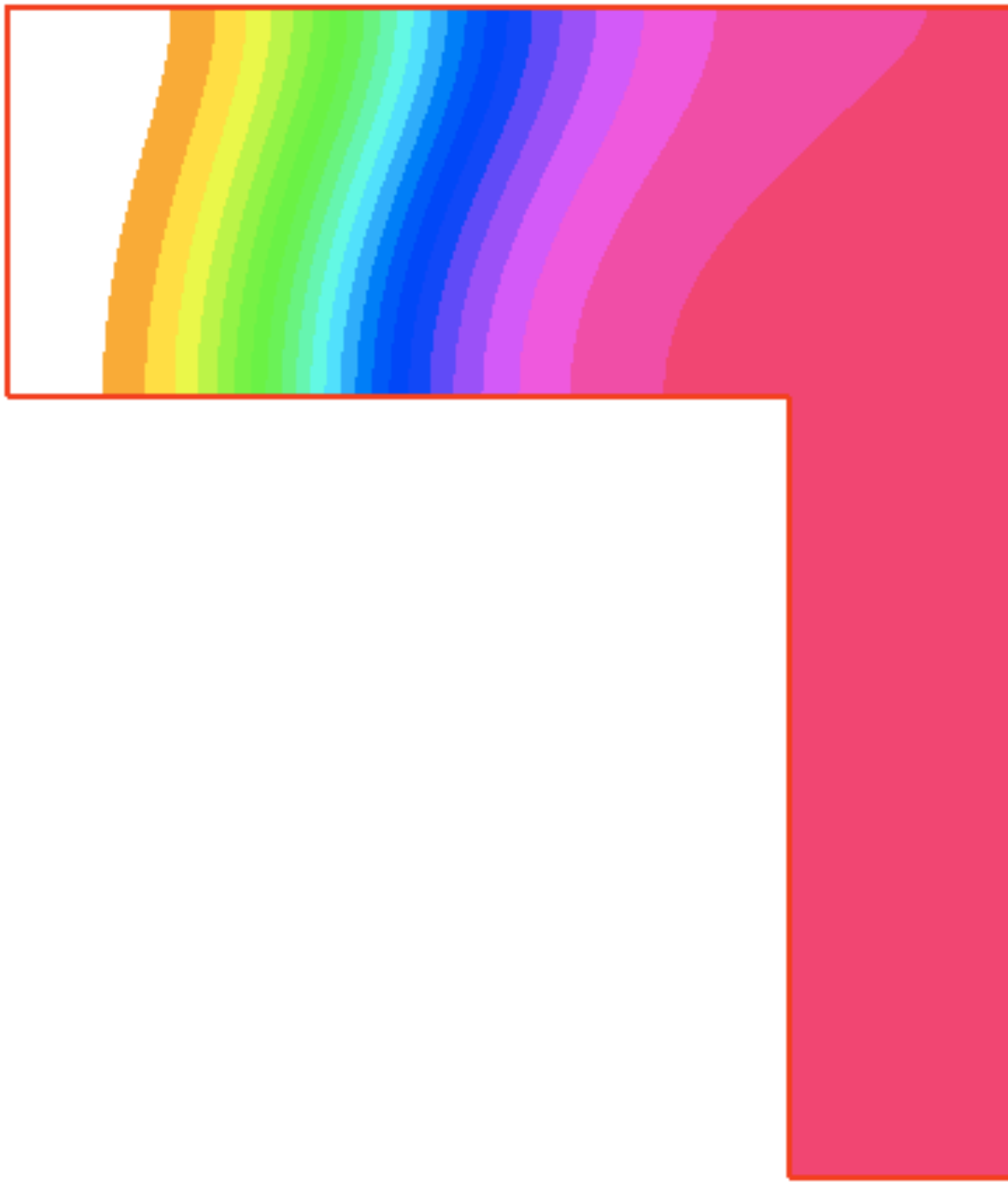}\\
\hspace{0.1cm} (a) & \hspace{0.1cm}(b) & \hspace{0.1cm}(c)
\end{tabular}}}
\caption{2D FKPP reference solution for the T-shape test case, at different times: (a) $t=0$, (b) $t=T_{max}/2$, (c) $t=T_{max}$.}
\label{fig::FKPP2D_Tref}
\end{figure}
\begin{figure}
\centerline{\hbox{\begin{tabular}{ccc}
\includegraphics[height=5.0cm,width=5.5cm]{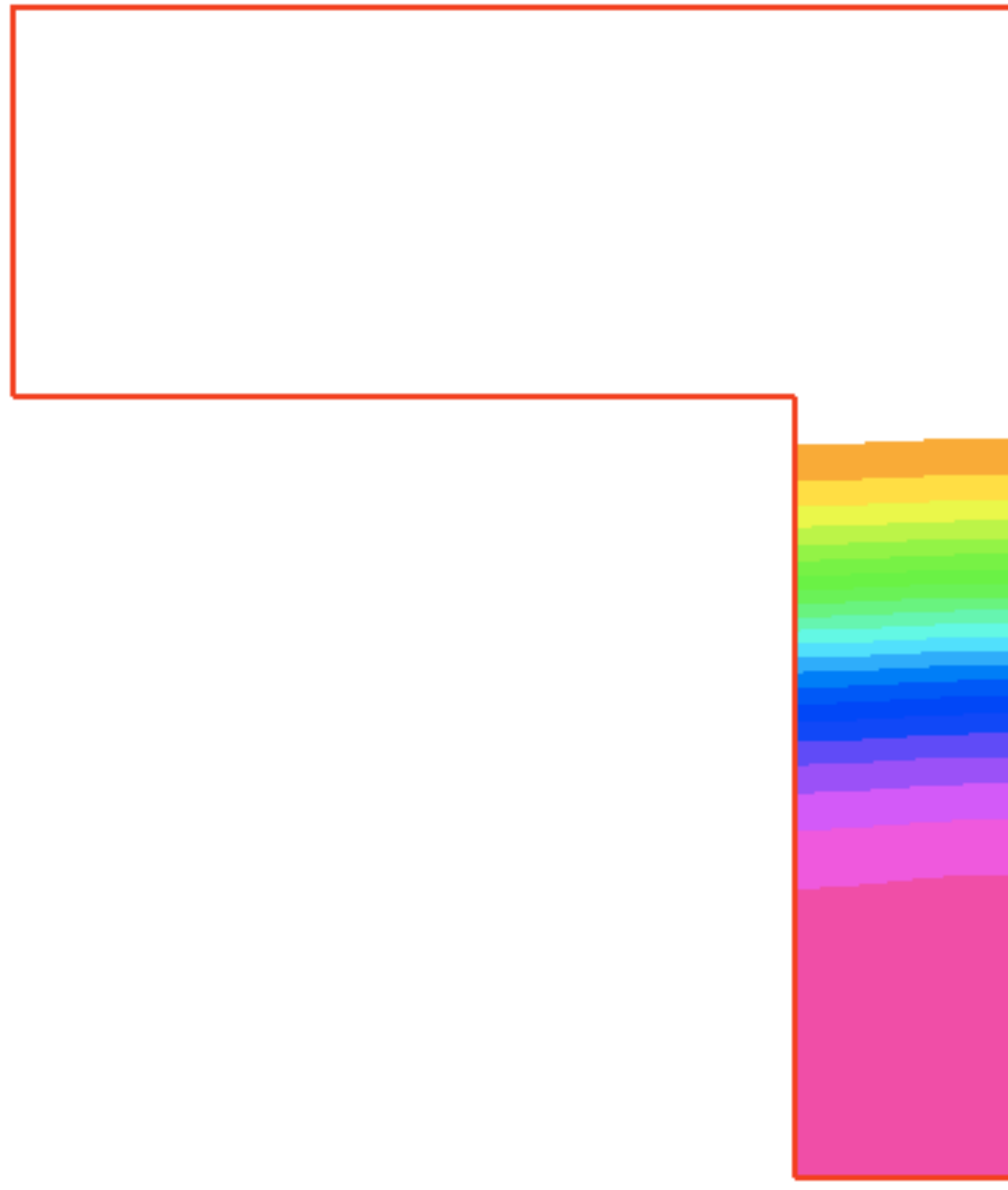} &
\includegraphics[height=5.0cm,width=5.5cm]{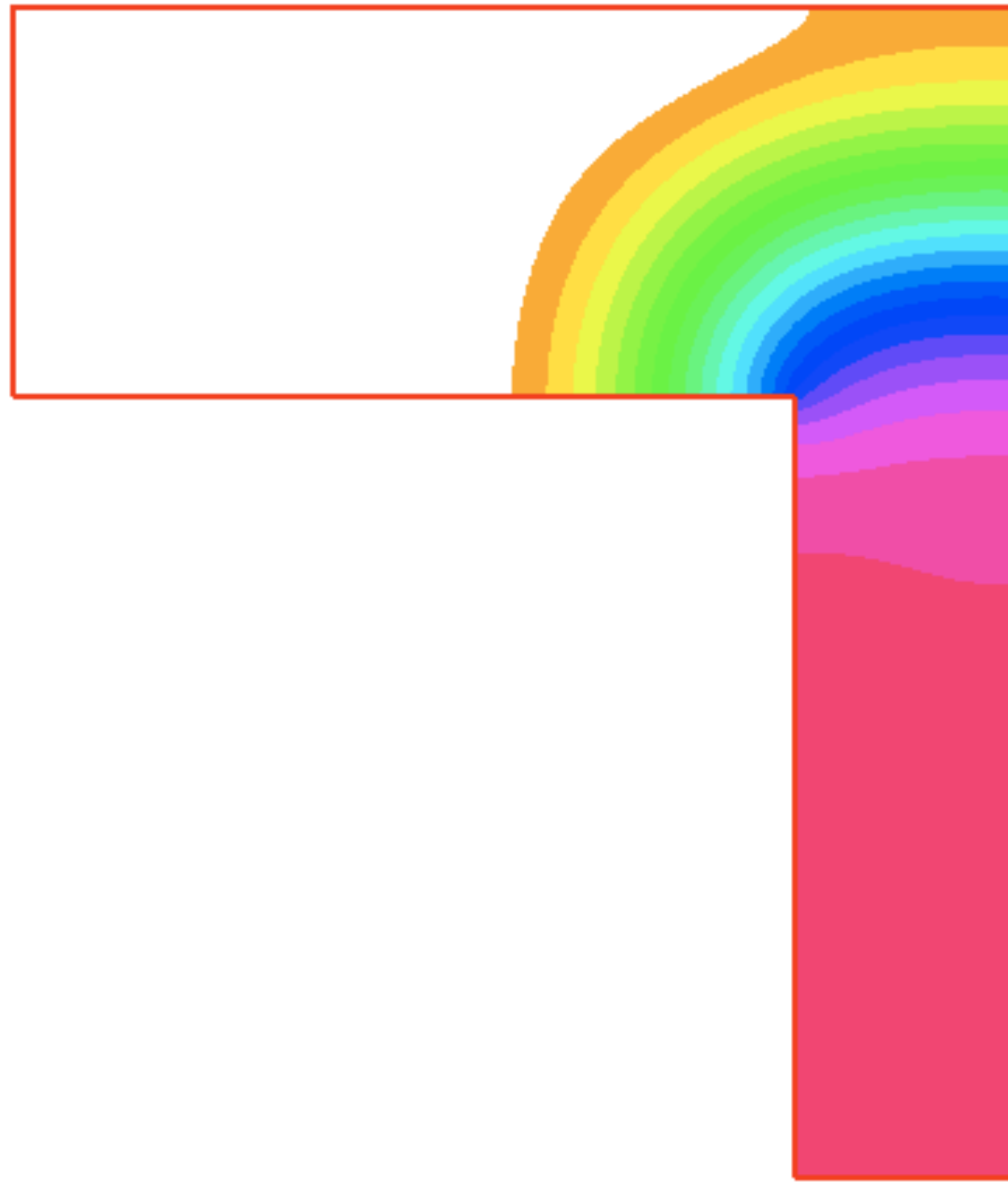} &
\includegraphics[height=5.0cm,width=5.5cm]{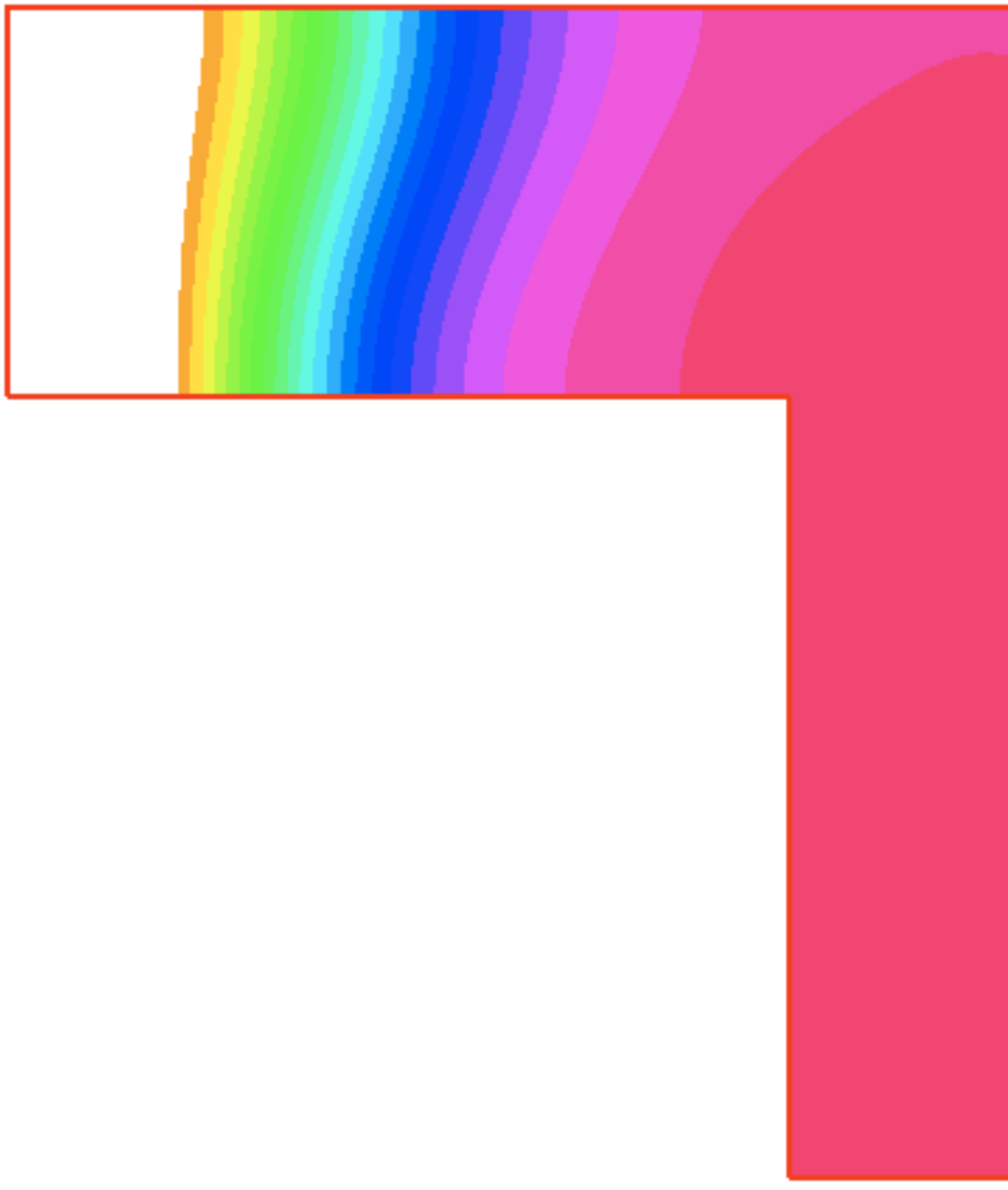}\\
\hspace{0.1cm} (a) & \hspace{0.1cm}(b) & \hspace{0.1cm}(c)
\end{tabular}}}
\caption{2D FKPP ROM solution for the T-shape test case, obtained with $N_M=30,\ \chi=1$ at different times: (a) $t=0$, (b) $t=T_{max}/2$, (c) $t=T_{max}$.}
\label{fig::FKPP2D_Talp}
\end{figure}
The results are shown for $N_M = 30$ modes. The $L^2$ relative error of the solution stays under $10\%$ for the whole simulation.
The average on the evolution simulated is $\overline{\varepsilon}_{L2} = 0.0279$. In Fig.\ref{fig::FKPP2D_Tref} three snapshots of the reference solution are shown at $t=0,T_{max}/2,T_{max}$. 
At the same time instants, the solution obtained by reconstructing the ROM solution is shown in Fig.\ref{fig::FKPP2D_Talp}. The dynamics is well recovered, the front position and shape are well rendered all along the evolution. We notice some inaccuracies concerning the front shape when the splitting occurs (see Fig.\ref{fig::FKPP2D_Talp}.(b)).  

\begin{rem}
Equation~\eqref{eq:beta-fkpp} is a vector logistic equation, whose stability of course depends on the respective influence of first order and quadratic terms.
The larger $\nu$ the larger is the number of modes needed to have $\nu<\lambda_p$. Depending on the problem symmetries and on the form of tensor $T$, we noticed that an unstable behavior could occur for  some subspaces of modes. We tested numerically the stability properties but a more careful analysis is in order. Roughly speaking, the larger $\nu$, the larger is the number of modes that have to be considered to have a stable integration of the ODE.
\end{rem}

\section{Comparison with the Semi Classical Signal Analysis }
\label{sec:reduced-approx-scsa}

In a preliminary version of this work~\cite{gerbeau:hal-00752810}, an approximation based on the Semi Classical Signal Analysis (SCSA) was used. Although it is less general than the approximation by the eigenmodes basis, it may be interesting in some applications. In this Section both approaches are compared.  

The SCSA was proposed in \cite{laleg-08}, analyzed in \cite{laleg-crepeau-sorine-12}, and successfully used for different applications to signal analysis in hemodynamics~\cite{laleg-medigue-sorine-07,laleg-medigue-van-de-louw-10}. It partially relies on the results by Lax and Levermore (see~\cite{Levermore_1} and~\cite{Levermore_2}). It consists in only keeping the eigenmodes of \eqref{eq:eigenfunction-def} corresponding to the negative eigenvalues $(\lambda_n)_{n=1... N_-}$ to approximate $u$ by the Deift-Trubowitz formula:
\begin{equation}
	\label{eq:u-deift-trubo}
	\tilde u({x}) = \chi^{-1}\sum_{m=1}^{N_-} \kappa_m \phi_m^2,
\end{equation}
with $\kappa_m = \sqrt{-\lambda_m}$. It clearly appears that this approach is limited to nonnegative signal\footnote{If $u({x})$ is not nonnegative, it is replaced by $u({x}) - \min_{x\in\Omega} u(x)$}. The parameter $\chi>0$ is chosen in order to reach the desired accuracy. For large values of $\chi>0$, the representation is more accurate, but also more expensive since the number of negative eigenvalues is larger. This decomposition is exact for a certain class of functions, called reflectionless potentials in physics. In the special case of the KdV equation, it corresponds to the decomposition of the solution in solitons. It has been shown in~\cite{laleg-medigue-sorine-07} that the artery blood pressure and flow rate can be accurately approximated with only a few modes with this formula. 

The approximation by the eigenfunctions and the SCSA are compared through their relative $L^2$ error
$\varepsilon^2_{L2} := \frac{\int_{\Omega} (u - \tilde{u})^2 \ d\Omega}{\int_{\Omega} u^2 \ d\Omega}$, where $u$ is the function that has to be approximated, and $\tilde{u}$ is obtained either by~\eqref{eq:u-modal} or \eqref{eq:u-deift-trubo}. 

\subsection{Static signals approximation}
The first tests deal with the approximation of given signals, without considering any dynamics. 

\paragraph{Realistic blood flow signal}
A first example is proposed on a realistic aortic flow. On this kind of signals, the SCSA \eqref{eq:u-deift-trubo} performs usually better than the approximation based on the eigenfunctions \eqref{eq:u-modal}.  

The parametric space $\chi : =[10^2,5 \ 10^3]$ was uniformly sampled. The maximum number of solitons (eigenfuctions squared) was, for each value of $\chi$, the number of negative eigenvalues. For the eigenfunction reconstruction, the approximation error was monitored up to $N_M = 40$ modes. 

\begin{table}
\begin{center}
\begin{tabular}{cccccc}
\hline
& $N_M$ & $\chi_s$  &$\varepsilon_{s}$ & $\chi_e$ & $\varepsilon_{e}$ \\
\hline
& $5$ & $1.04\ 10^3$ & $6.37 \ 10^{-2}$ & $1.38\ 10^3$ & $7.83\ 10^{-2}$\\
& $6$ & $1.63\ 10^3$ & $4.94 \ 10^{-2}$ & $1.88\ 10^3$ & $6.79\ 10^{-2}$\\ 
& $7$ & $2.28\ 10^3$ & $3.99 \ 10^{-2}$ & $3.07\ 10^3$ & $5.38\ 10^{-2}$\\
& $8$ & $3.17\ 10^3$ & $3.43 \ 10^{-2}$ & $3.81\ 10^3$ & $4.73\ 10^{-2}$\\ 
& $9$ & $4.21\ 10^3$ & $3.16 \ 10^{-2}$ & $3.76\ 10^3$ & $4.35\ 10^{-2}$\\  
& $10$ & $4.31\ 10^3$ & $2.82 \ 10^{-2}$ & $3.61\ 10^3$ & $4.30\ 10^{-2}$\\ 
\hline
\end{tabular}
\caption{Errors of the two representations as a function of the number of modes used ($N_M$). The columns $\chi_{s}$ and $\chi_e$ are the values of the parameter for which the error of the soliton and eigenfunction reconstructions is the smallest one, $\varepsilon_{s}$ and $\varepsilon_{e}$ are the errors.}
\label{table::staticErrors}
\end{center}
\end{table}
In Table \ref{table::staticErrors} the errors for the two reconstructions are reported. In particular, for a fixed $N_M$ the optimal $\chi$ and the associated error are written. The two representations give similar results. However, the soliton reconstruction is slightly better and, for certain values of the parameter $\chi$, we need to increase the number of modes for the eigenfunction reconstruction in order to have the same performances as with the solitons. 
\begin{figure}
\centerline{\hbox{\begin{tabular}{cc}
\includegraphics[height=7.5cm]{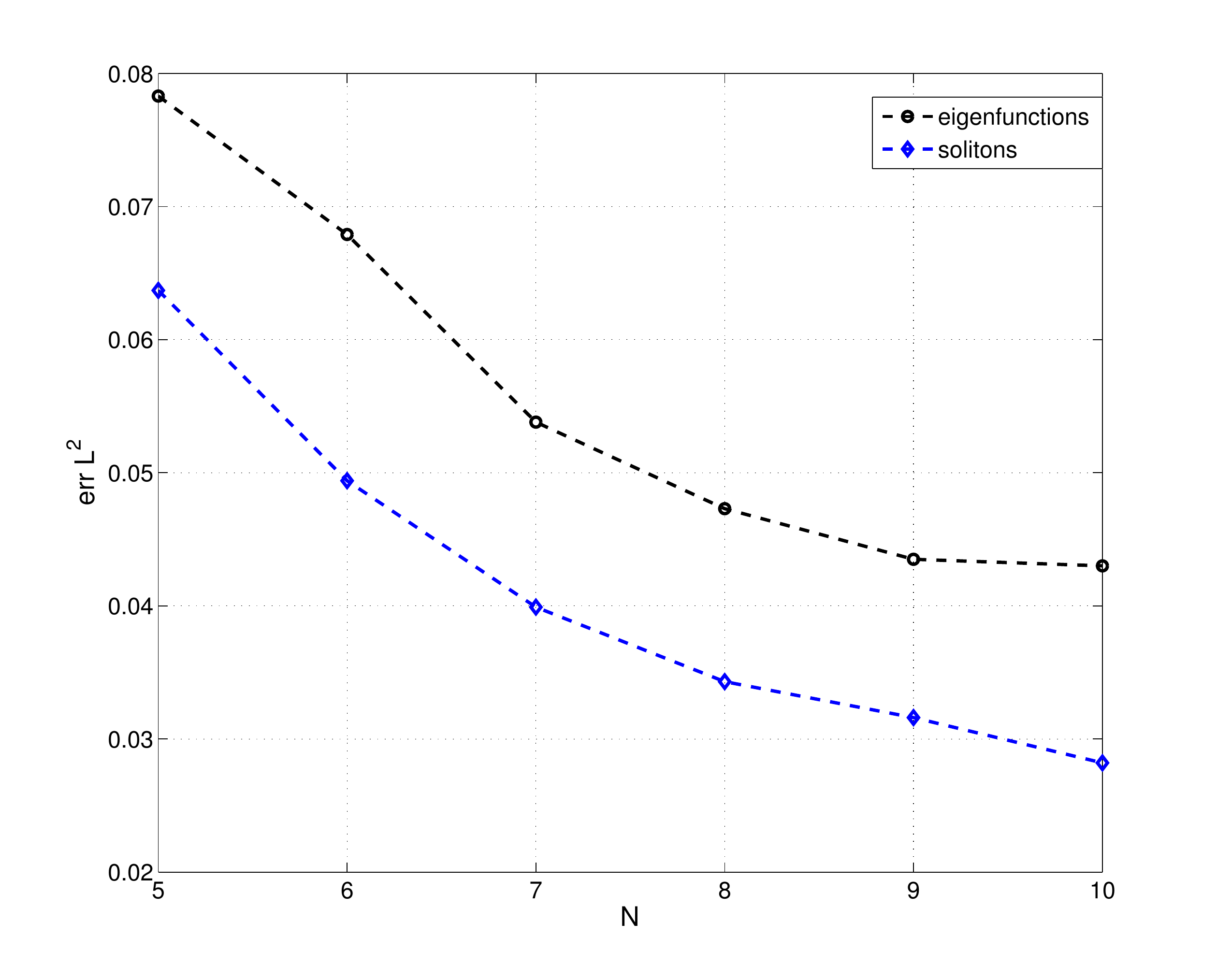} \\ (a) \\
\includegraphics[height=7.5cm]{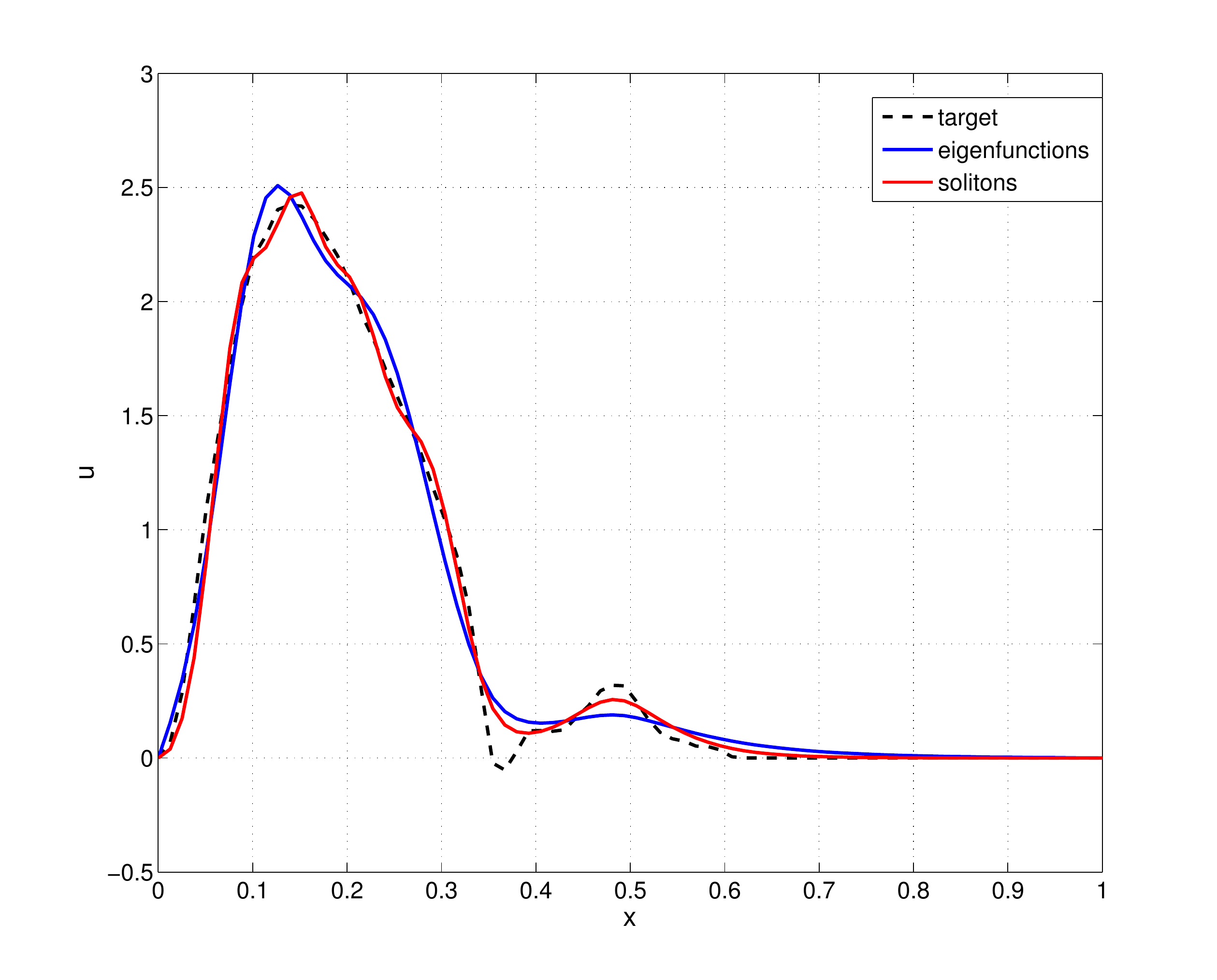}\\ (b) 
\end{tabular}}}
\caption{ (a) Errors in $L^2$ norm for optimal $\chi$, varying the number of modes, see Table\ref{table::staticErrors} (b) Comparison : in blue the eigenfunctions reconstruction, in red the eigenfunctions squared one, black-dashed is the target solution. $N=5,\chi=1000$. }
\label{fig::cardiacFlux}
\end{figure} 
In Fig.\ref{fig::cardiacFlux}.(a) the errors in $L^2$ norm are shown as a function of the number of modes, for an optimal choice of the parameter $\chi$ (see Table \ref{table::staticErrors}). In Fig.\ref{fig::cardiacFlux}.(b) the reconstructions are compared: the dot-dashed line, in black, is the target, the reconstruction based on the eigenfunction expansion is plotted in blue, the soliton one in red. The two reconstructions are similar but, for a given number of modes, the one based on solitons is better at capturing the features of the signal. 

\begin{figure}
\centerline{\hbox{\begin{tabular}{cc}
\includegraphics[height=7.0cm]{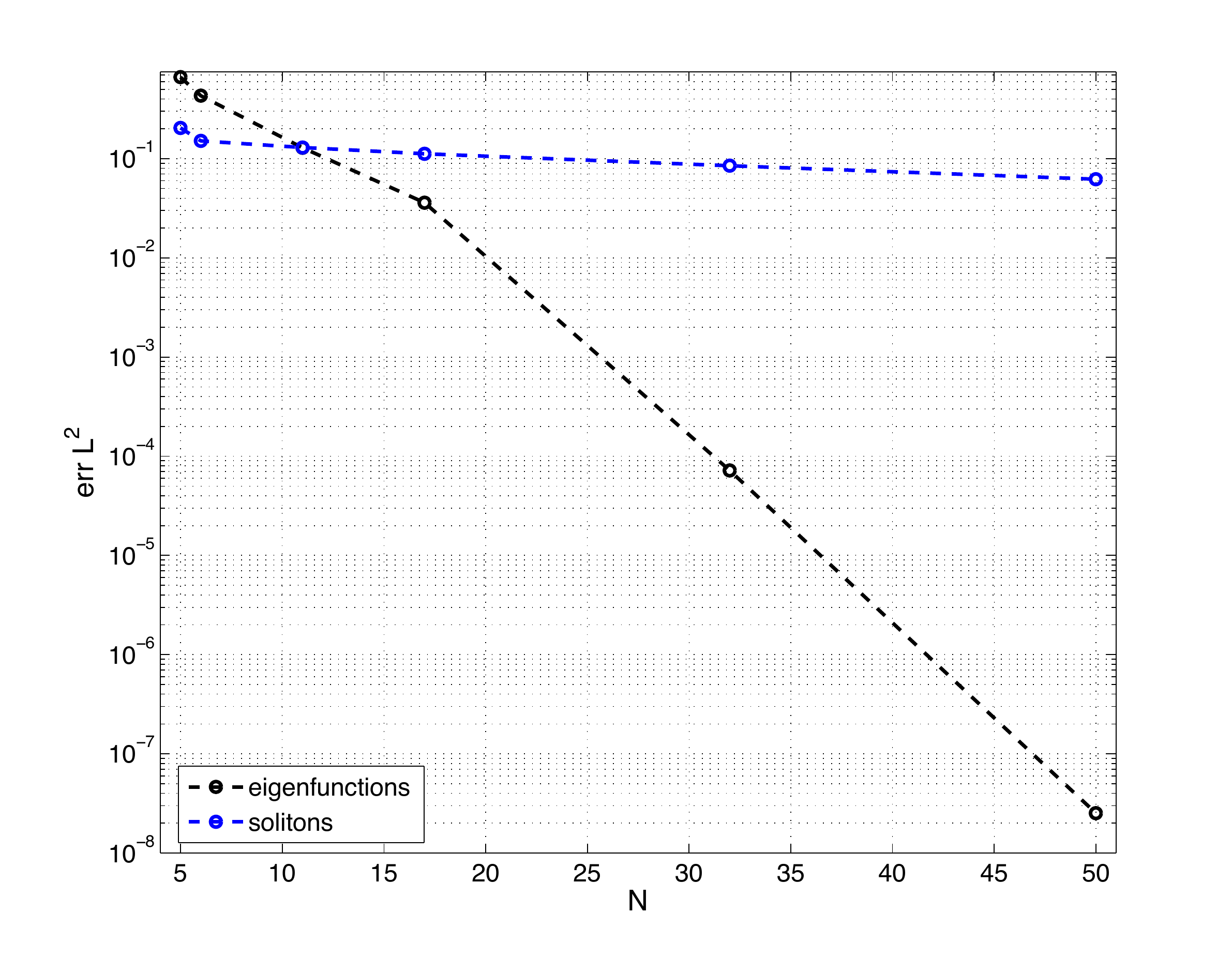} \\ (a) \\
\includegraphics[height=7.0cm]{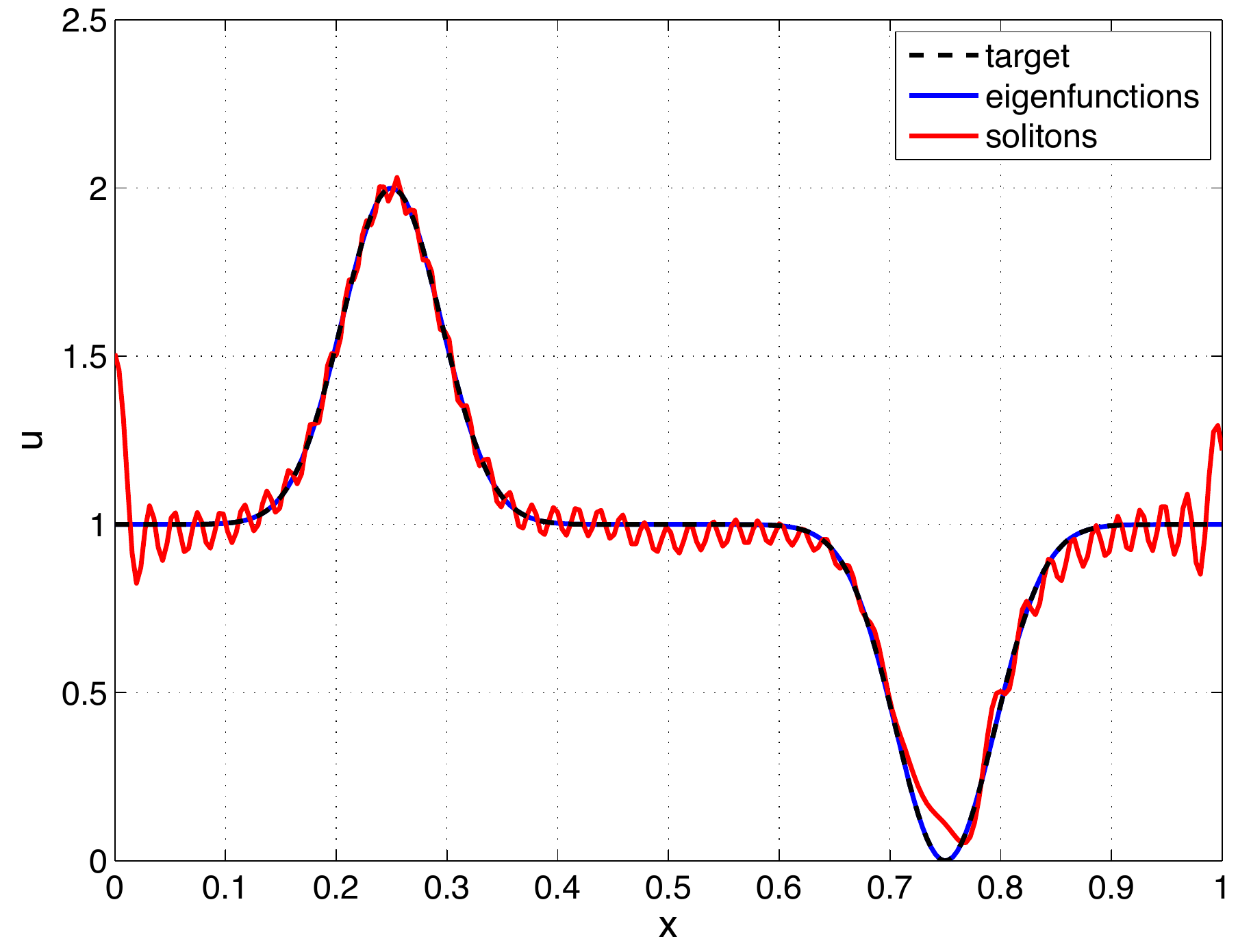}\\ (b) 
\end{tabular}}}
\caption{(a) Errors in $L^2$ norm varying the number of modes in semi-logarithmic scale (b) Comparison : in blue the eigenfunctions reconstruction ($\chi = 250$, $N=50$), in red the eigenfunctions squared one ($\chi = 2.5e4$, $N=50$), black-dashed is the target solution. }
\label{fig::doubleGauss}
\end{figure}

\paragraph{Double gaussian profile}
We consider a target function defined by $u(x) = \exp(-250(x-0.25)^2) - \exp(-250(x-0.75)^2)$. Note that it has a negative part, the approximation based on solitons cannot be used in that case. The function was therefore translated to be nonnegative and then both the reconstructions were tested (see Figure~\ref{fig::doubleGauss}).
In Fig.\ref{fig::doubleGauss}.(a) the error in $L^2$ norm is shown in semi-logarithmic scale as a function of the number of modes used. The eigenfunction reconstruction (in black), built by setting $\chi=250$, converges very fast, while the reconstruction based on the solitons (in blue) converges poorly. In Fig.\ref{fig::doubleGauss}.(b) a comparison of the reconstructions is shown when $N=50$, that confirms that the eigenfunctions squared reconstruction is not well adapted in this case.

In conclusion, the approximation by eigenfunctions has a clear advantage of generality. Nevertheless, the approximation by solitons may be interesting for some specific signals and further studies would be useful to better understand its properties.

\subsection{KdV equation}
\label{sec:kdv-SCSA}
The KdV equation, when the solution is an $n$-solitons, is a typical example in which the expansion of the solution as sum of eigenfunctions squared performs better, the error being merely due to space discretization. 
Indeed, in this case, as well as for other integrable systems, the proposed approach is a numerical discretization of Lax pairs, whose analytical expression for KdV (see for instance \cite{Ablowitz_1}) reads:
\begin{eqnarray}
& \mathcal{L}(u) \cdot = -\partial^2_x \cdot -u  \cdot  \ , \label{KdV_Operators-L}\\
& \mathcal{M}(u)  \cdot  = 4 \partial^3_x \cdot+3u\partial_x \cdot + 3\partial_x(u  \cdot  ). \label{KdV_Operators-M}
\end{eqnarray}
Doing as if the Lax pair was unknown, we consider the one-soliton and the three-soliton propagations.
Let us denote the number of negative eigenvalues of the Schr\"{o}dinger operator by $N_{-}$. The soliton reconstruction $u = \sum_{i=1}^{N_{-}}\alpha_i \phi_i^2$ is assumed and the ALP algorithm  rederived accordingly, only few changes being necessary.

The soliton expansion is injected into the Eq.\eqref{eq:KdVAnalytical}, leading to:
\begin{equation}
\sum_{i=1}^{N_{-}} \partial_t \alpha_i \phi_i^2 + \alpha_i \left(\partial_t (\phi_i^2) - 4\lambda_i \partial_x (\phi_i^2)\right) + 4 \sum_{i,j=1}^{N_{-}} \alpha_i\alpha_j \left( \phi_i^2 \partial_x(\phi_j^2) - \phi_j^2 \partial_x(\phi_i^2)\right) = 0. 
\end{equation}
Since the last term is a quadratic symmetric form of a skew-symmetric term, it is equal to zero and the equation reduces to:
\begin{equation}
\sum_{i=1}^{N_{-}} \partial_t \alpha_i \phi_i^2 + \alpha_i \left(\partial_t (\phi_i^2) - 4\lambda_i \partial_x (\phi_i^2)\right) = 0,
\end{equation}
which highlights some properties of the solution. By projecting this equation on the basis, an evolution ODE for the coefficients is obtained:
\begin{equation}
\partial_t\alpha_i + 2 \sum_{j=1}^{N_{-}} (M_{ij} - 4 \lambda_j D_{ij}) \alpha_j=0,
\label{eq:alpha-KdV}
\end{equation}
where $D_{ij} = \langle \partial_x\phi_j, \phi_i \rangle$. The evolution of the matrix $D$ is governed by 
\begin{equation}
	\label{eq:D-dyn-KdV}
\dot D + [D,M] = 0.
\end{equation}
The basis evolution is accounted for by using the Eq.$\eqref{eq:rom-dynamics}_{2-4}$ and, in this case, $\gamma = \gamma(\alpha) = 8 \sum_{j=1}^{N_{-}}\lambda_j \alpha_j D_{ij}$ is substituted to Eq.$\eqref{eq:rom-dynamics}_{5}$.

\paragraph{One-soliton solution}
As $u_0$ is the initial datum of the one-soliton propagation and $\chi=1$ provides the analytical expression for $\mathcal{L}(u)$ in the case of the KdV equation (see Eq.(\ref{KdV_Operators-L})), only one eigenvalue belongs to the discrete spectrum and the corresponding mode squared is exactly $u_0$, up to discretization errors ($10^{-4}$ in $L^2$ norm for the present case). 
\begin{figure}
\centerline{\hbox{\begin{tabular}{cc}
\includegraphics[height=6.5cm]{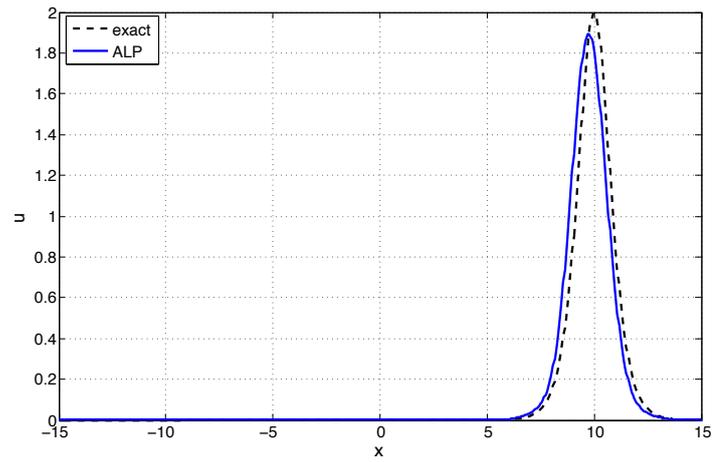} \\ (a) \\
\includegraphics[height=6.5cm]{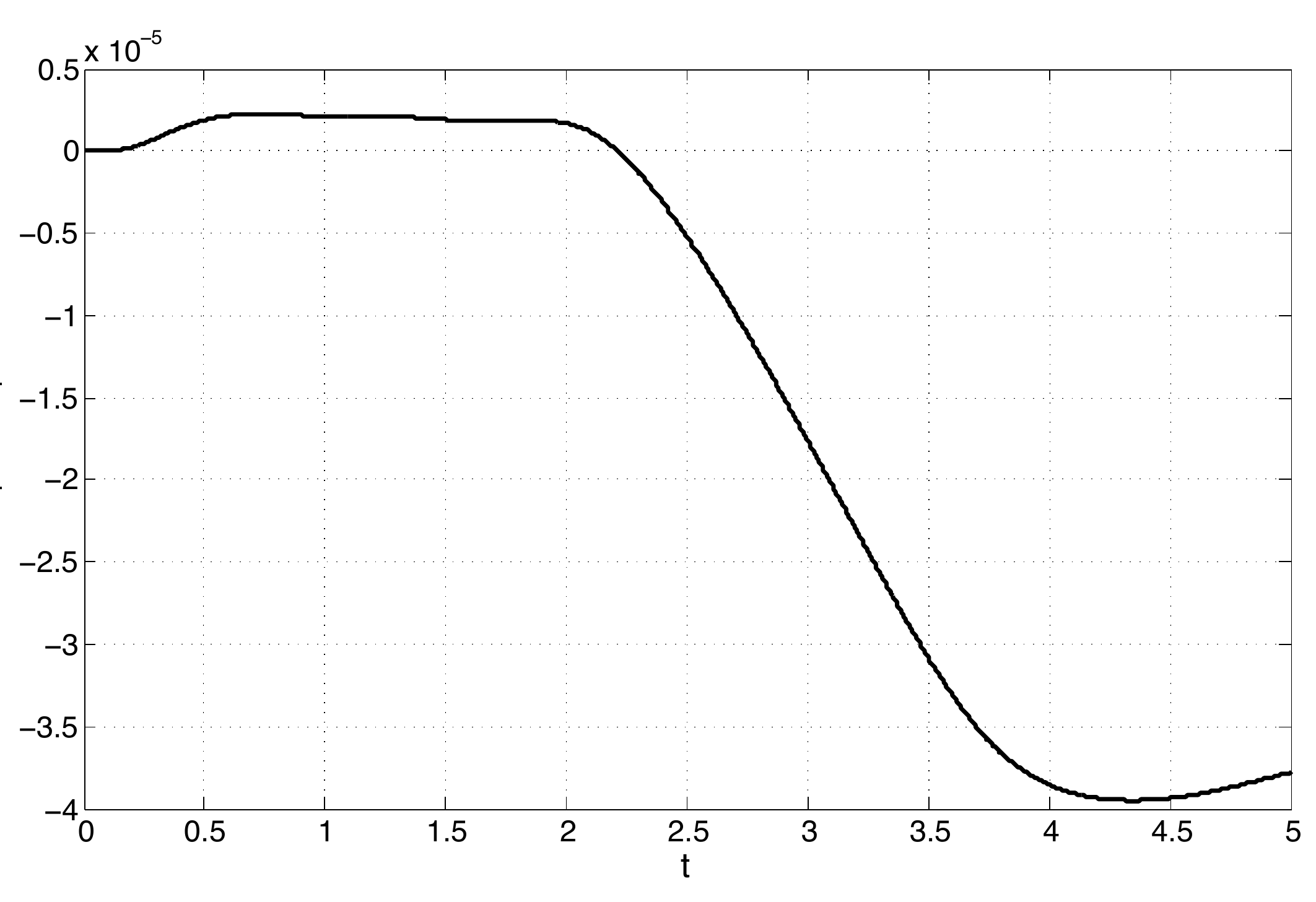}\\ (b)
\end{tabular}}}
\caption{(a) Comparison between the exact solution and the solution obtained by integrating ALP with $N_M=25$, $N_{-}=1$, (b) error in time for $\alpha_1$ when ALP is integrated with $N_M=5$.}
\label{fig::KdV1}
\end{figure}

The ALP-ROM was integrated by using a $\delta t = T_{max}/2500$, varying the number of modes $N_M$ used to represent the operators. 

For the KdV equation, an analytical results holds for the coefficients: $\alpha_1(t) = \alpha_1|_{t=0}$.
The Reduced Order Model allows to recover this result: in Fig.\ref{fig::KdV1}.(b) the error between the exact and the simulated value for $\alpha_1$ is shown as function of time when only $N_M=5$ modes were used.
The error on $\alpha_1$ in the reduced space weakly depends upon $N_M$. 
The number of modes used to discretize the operators has an influence in the postprocessing stage, so that it affects the error between the analytical and the reconstructed solution.  

\begin{table}
\begin{center}
\begin{tabular}{cccc}
\hline
& $N_M$ & $\overline{\varepsilon_{L2}}$ & $\max_t(\varepsilon_{L2})$ \\
\hline
& $26$ & $0.0817$ & $0.1774$ \\
& $28$ & $0.0524$ & $0.1141$ \\
& $30$ & $0.0374$ & $0.0826$ \\
& $32$ & $0.0302$ & $0.0690$ \\
& $34$ & $0.0269$ & $0.0637$ \\
& $36$ & $0.0184$ & $0.0517$ \\
\hline
\end{tabular}
\caption{Error indicators for the KdV one-soliton test case as a function of the number of modes used to discretize the equations: first column $N_M$ is the number of modes, the second and the third one the time average and the maximum of the $L^2$ error.}
\label{table::KdV1-solit}
\end{center}
\end{table}
In Table \ref{table::KdV1-solit}, the error indicators for this case are shown as function of the number of modes used to discretize the operators. The qualitative agreement between the reconstructed solution and the analytical one is shown, at final time, in Fig.\ref{fig::KdV1}.(a), for $N_M=25$: all the features of the wave are well captured by the reduced order solution.

\paragraph{Three-soliton solution}

The spectral problem is solved at initial time and, by setting $\chi=1$, three distinct eigenvalues are found in the negative part of the spectrum.
This is in agreement with the analytical results and highlights the ability to decompose a traveling (non-linearly interacting) waves system in its basic components, and propagate them separately.

\begin{table}
\begin{center}
\begin{tabular}{cccc}
\hline
& $N_M$ & $\overline{\varepsilon_{L2}}$ & $\max_t(\varepsilon_{L2})$ \\
\hline
& $28$ & $0.0237$ & $0.0332$ \\
& $32$ & $0.0123$ & $0.0204$ \\
& $36$ & $0.0098$ & $0.0185$ \\
& $40$ & $0.0047$ & $0.0102$ \\
& $44$ & $0.0023$ & $0.0086$ \\
& $48$ & $0.0015$ & $0.0038$ \\
\hline
\end{tabular}
\caption{Error indicators for the KdV three-soliton test case as a function of the number of modes used to discretize the equations: first column $N_M$ is the number of modes, the second and the third one the time average and the maximum of the $L^2$ error.}\label{table::KdV3}
\end{center}
\end{table}
The results are similar to those obtained for the simpler one-soliton case. 
In particular, the coefficients $\alpha_{1,2,3}$ do not vary in time up to $10^{-4}$, so that the error in the reduced space is negligible and the analytical result is recovered.
The error in the high dimensional space is governed by the number of modes $N_M$ used for the discretization of the operators and in the post-processing stage. The errors are shown in Table \ref{table::KdV3}. 

\section{Conclusions and perspectives}

We have proposed a new reduced-order model technique, called ALP, consisting of three stages. First, a set of orthonormal eigenfunctions of a linear Schr\"odinger operator associated with the initial condition is computed. Second, a projection of the PDE on the time dependent basis of the reduced order space is solved. Third, the solution is reconstructed on the full order space by propagating the reduced order basis in time with an approximation of a Lax operator. Interestingly, it is not necessary to perform the reconstruction stage to solve the equation in the reduced order space. 

The method was successfully tested on the linear advection, the KdV and the FKPP equations in 1D and 2D. It seems to be well-adapted to systems modeling propagation phenomena. Unlike other reduced-order methods, it does not rely on an off-line computation of a large data set of solutions.

The application of ALP to other problems is currently under investigation, in particular to a set of Euler equations modeling a network of arteries and to cardiac electrophysiology problems. Many questions would deserve further investigations: the number of modes could be adapted along the resolution, for example based on the indicator \eqref{eq::frobErr}; other operators than the Laplacian might used for operator $\mathcal{L}$; other time schemes could be used to solve the reduced order dynamics \eqref{eq:rom-dynamics} or the modes propagation \eqref{eq:phi-with-residual}; a more precise reconstruction method could be devised; the role of parameter $\chi$ should be further investigated; the scheme could be extended to handle non-polynomial nonlinearity; etc. This will be the subject of future works.

\section*{References}


\clearpage
\end{document}